%% file: longsurvey.tex
\documentclass[11pt]{article}
\usepackage{diagrams}
\input{format} \input{mathdefs} \input{thmdefs}

\newtheorem{probl}{Problem} 
\newenvironment{problem}{\begin{probl} \rm}{\end{probl}}
\newcommand{\Min}{\operatorname{Min}}
\newcommand{\MIN}{\operatorname{MIN}}
\newcommand{\COL}{\operatorname{COL}}
\newcommand{\Max}{\operatorname{Max}}
\newcommand{\MAX}{\operatorname{MAX}}
\newcommand{\ROW}{\operatorname{ROW}}
\newcommand{\SB}{{\cal SB}}

\newcommand{\CB}{{\cal CB}}
\newcommand{\cb}{{\mathrm{cb}}}
\newcommand{\PM}{\operatorname{PM}}
\begin{document}
\title{Applications of operator spaces to abstract harmonic analysis}
\author{\it Volker Runde}
\date{}
\maketitle
\begin{abstract}
We give a survey of how the relatively young theory of operator spaces
has led to a deeper understanding of the Fourier algebra of a locally
compact group (and of related algebras).
\end{abstract}
\section*{Introduction}
Abstract harmonic analysis is the mathematical discipline concerned with the study of locally compact groups and of the spaces and algebras associated with them.
\par
The framework of abstract harmonic analysis was ready when A.\ Weil proved the existence (and uniqueness) of (left) Haar measure on an arbitrary locally compact group $G$ (\cite{Wei}) after A.\ Haar had previously
dealt with the case where $G$ was also supposed to be separable and metrizable (\cite{Haa}). With the existence of Haar measure on $G$, one can then, of course, consider the $L^p$-spaces $L^p(G)$ for $p \in [1,\infty]$
and, in particular, equip $L^1(G)$ with a convolution product turning it into a Banach algebra. The group algebra $L^1(G)$ is a complete invariant for $G$ in the sense that, if $H$ is another locally compact group
such that the Banach algebras $L^1(G)$ and $L^1(H)$ are isometrically isomorphic, then $G$ and $H$ are topologically isomorphic (\cite{Wen}). Hence, the powerful theory of Banach algebras can be applied to
study locally groups; the first monograph to treat abstract harmonic analysis in a Banach algebraic context was \cite{Loo}. 
\par
A question related to the existence of Haar measure is whether, for a particular $G$, there is a non-zero, positive, linear functional on $L^\infty(G)$ that is invariant under left translation: this question was first investigated by 
J.\ von Neumann (\cite{Neu}). He called (discrete) groups for which such functionals exists ``Gruppen von endlichem Ma{\ss}''. Nowadays such groups are called {\it amenable\/} following M.\ M.\ Day (\cite{Day}). Both compact and abelian 
groups are amenable whereas the free group in two generators isn't. The theory of amenable, locally compact groups is expounded in \cite{Gre}, \cite{Pie}, and \cite{Pat}.
\par
In his seminal memoir \cite{Joh1}, B.\ E.\ Johnson characterized the amenable locally compact groups through a cohomological vanishing condition of their $L^1$-algebras. This triviality condition makes sense for every
Banach algebra and is used to characterize the (rich) class of amenable Banach algebras. Since its introduction by Johnson, the concept of an amenable Banach algebra has turned out to be fundamental. For $\cstar$-algebras, for example,
amenability is equivalent to the pivotal notion of nuclearity (see \cite[Chapter 6]{LoA}). At about the same time Johnson defined amenable Banach algebras, A.\ Ya.\ Helemski\u{\i} in Moscow began to systematically develop the subject of 
topological homology, i.e.\ of homological algebra with functional analytic overtones added (see \cite{Hel} for an account). It turns out that amenability for Banach algebras fits nicely into this framework.
\par
If $G$ is abelian with dual group $\hat{G}$, the Fourier ($=$ Gelfand) transform maps $L^1(G)$ onto a subalgebra of ${\cal C}_0(\hat{G})$ which is denoted by $A(\hat{G})$. It follows from Plancherel's theorem that $A(\hat{G})$ consists
precisely of those functions on $\hat{G}$ which are the convolution product of two $L^2$-functions: this was used by P.\ Eymard to define the {\it Fourier algebra\/} $A(G)$ for arbitrary $G$ (\cite{Eym}). The first to characterize
properties of $G$ in terms of the Banach algebra $A(G)$ was H.\ Leptin, who proved that $G$ is amenable if and only if $A(G)$ has a bounded approximate identity (\cite{Lep}). The tempting conjecture, however, that $A(G)$ is amenable
(as a Banach algebra) if and only if $G$ is amenable (as a group) turns out to be false --- in \cite{JohFourier}, Johnson showed that $A(G)$ fails to be amenable for certain compact $G$.
\par
In 1995, Z.-J.\ Ruan published a result that would shed new light on $A(G)$ and initiated a completely novel approach to studying the Fourier algebra (\cite{RuaAG}). The key to this approach is the still fairly young theory of (abstract) 
operator spaces. Originally, an operator space was defined to be a closed subspace of ${\cal B}(\Hilbert)$ for some Hilbert space $\Hilbert$ (see, e.g., \cite{Pau1}). In his ground breaking paper \cite{Rua}, Ruan characterized operator spaces
by means of two simple axioms that involve norms on all spaces of matrices over the given space. The advantages of this axiomatic approach over the old, concrete one are manifold: for instance, it allows for the development of a duality
theory that parallels the duality of Banach spaces in many aspects (see, e.g., \cite{ER}). In particular, the (Banach space) dual of an operator space is again an operator space in a canonical manner. The Fourier algebra $A(G)$ can be 
canonically identified with the unique predual of the group von Neumann algebra $\VN(G)$ and thus carries a natural operator space structure. The notion of an amenable Banach algebra easily adapts to the operator space context and yields
what is called operator amenability (\cite{RuaAG}). It turns out that this operator space theoretic variant of amenability for Banach algebras is the ``right'' one when it comes to dealing with Fourier algebras: $A(G)$ is amenable if
and only if $G$ is amenable (\cite{RuaAG}).
\par
Since the publication of \cite{RuaAG}, various authors have successfully used the operator space structure of $A(G)$ to match homological properties of that algebra with properties of $G$ (\cite{Ari}, \cite{ARS}, \cite{Spr}, and 
\cite{Woo}). Even if one is interested in $A(G)$ only as a Banach algebra, operator space techniques turn out to be valuable. In \cite{FKLS}, B.\ E.\ Forrest et al.\ use operator space methods to characterize, for amenable $G$, those closed ideals of $A(G)$ which have a bounded approximate identity, and in \cite{FR}, Forrest and the author  characterize those locally compact groups $G$ for which $A(G)$ is amenable, also relying on
operator space techniques.
\par
The present article is intended as a survey of operator space methods in the investigation of $A(G)$ (and related algebras). We suppose that the reader is fluent in basic functional analysis (including the fundamentals of operator algebras)
and is willing to look up some of the background from abstract harmonic analysis and operator spaces in the references given. We often omit proofs altogether or present them rather sketchily; the word ``proof'' can therefore often just
mean ``idea of a proof''. There is some overlap with the recent, much less detailed article \cite{RunSurv}, which we have striven to keep to a minimum.
\section{Locally compact groups and their group algebras} \label{aha}
A locally compact group is a group $G$ equipped with a locally compact Hausdorff topology such that the maps
\[
  G \times G \to G, \quad (x,y) \mapsto xy \qquad\text{and}\qquad G \mapsto G, \quad x \mapsto x^{-1}
\]
are continuous.
\par
Trivially, every group equipped with the discrete topology is locally compact. Also, every Lie group is locally compact. This latter class immediately supplies us with a multitude of (non-discrete) examples: $\reals^N$, $\torus^N$ where 
$\torus := \{ c\in \comps : |c| = 1 \}$, and matrix groups such as $\operatorname{GL}(N,\comps)$,  $\operatorname{SL}(N,\reals)$, etc.
\par
The probably most surprising fact about objects as general as locally compact groups is that there are still substantial theorems to prove about them. The starting point of abstract harmonic analysis (as opposed to classical harmonic
analysis) is the following theorem that was proved in its full generality by A.\ Weil (\cite{Wei}):
\begin{theorem} 
Let $G$ be a locally compact group. Then there is a non-zero, regular (positive) Borel measure $G$ --- {\rm (left) Haar measure} --- which is left invariant, i.e.\ $x B$ and $B$ have the same measure for all $x \in G$ and all Borel subsets 
$B$ of $G$, and unique up to a multiplicative, positive constant.
\end{theorem}
\par
If $G$ is a locally compact group and $B \subset G$ is a Borel set, we write $|B|$ for the Haar measure of $B$; integration with respect to Haar measure is denoted by $dx$.
\par
For discrete $G$, Haar measure is just counting measure; for $G = \reals^N$, it is $N$-dimensional Lebesgue measure; for $G = \torus$, it is arclength measure.
\par
It is not true that Haar measure is always right invariant, but there is a unique continuous group homomorphism $\Delta \!: G \to (0,\infty)$ such that $|Bx| = \Delta(x) |B|$ for each $x \in G$ and for each Borel set $B \subset G$.
For many groups, $\Delta \equiv 1$ holds, e.g., if $G$ is abelian (trivially), discrete (because Haar measure is counting measure), or compact (because $\Delta(G)$ must be a compact subgroup of $(0,\infty)$ and thus equal $\{ 1 \}$),
even though this is false for general $G$. 
\par
We now bring Banach algebras into the picture: For any $f, g\in L^1(G)$ define their {\it convolution product\/} $f \ast g \in L^1(G)$ by letting
\begin{equation} \label{conv}
  (f \ast g)(x) := \int_G f(y)g(y^{-1}x) \, dy \qquad (x \in G).
\end{equation}
This formula has to be read, of course, with the basic precautions: If $f$ and $g$ are $L^1$-functions on $G$, then the integral on the right hand side of (\ref{conv}) exists for almost all $x \in G$, only depends on the equivalence classes
of $f$ and $g$, respectively, and defines (almost everywhere) an $L^1$-function on $G$ denoted by $f \ast g$ --- all this follows easily from the Fubini--Tonelli theorem. It is routinely checked that the convolution product turns $L^1(G)$ into
a Banach algebra.
\par
Already in the introduction, we quoted the following theorem by J.\ G.\ Wendel (\cite{Wen}):
\begin{theorem} \label{wendel1}
Let $G$ and $H$ be locally compact groups. Then $L^1(G)$ and $L^1(H)$ are isometrically isomorphic if and only if $G$ and $H$ are topologically isomorphic.
\end{theorem}
\par
Consequently, every property of $G$ that can be expressed in terms of the locally compact group structure can be expressed in terms of $L^1(G)$. It is easy to see that $L^1(G)$ is commutative if and only if $L^1(G)$ is abelian and that
$L^1(G)$ has an identity if and only if $G$ is discrete. Other --- much less easily seen --- correspondences of this kind will be discussed later in this article.
\par
Even though $L^1(G)$ for non-discrete $G$ lacks an identity it has something almost as good:
\begin{theorem} \label{BAI}
Let $G$ be a locally compact group and let $\mathfrak U$ be a basis of neighborhoods of the identity of $G$. Furthermore, for each $U \in \mathfrak U$, let $e_U \in L^1(G)$ be positive with $\| e_U \|_1 = 1$ such that $\supp \, e_U \subset U$.
Then $( e_U )_{U \in {\mathfrak U}}$ is a\/ {\rm bounded approximate identity\/} for $L^1(G)$, i.e.\
\[
  f \ast e_U \to f \quad\text{and}\quad e_U \ast f \to f \qquad (f \in L^1(G)).
\]
\end{theorem}
\par
Moreover, the $L^1$-algebra of a locally compact group carries a natural involution:
\[
  f^\ast(x) := \frac{1}{\Delta(x)} \overline{f(x^{-1})} \qquad (f \in L^1(G), \, x \in G).
\]
It is easily checked that this involution is isometric.
\par
By a representation of $G$ on a Hilbert space $\Hilbert$, we mean a homomorphism from $G$ into the unitaries on $\Hilbert$ which is continuous with respect to the given topology on $G$ and the strong operator topology on ${\cal B}(\Hilbert)$.
Given any such representation $\pi \!: G \to {\cal B}(\Hilbert)$, we obtain a $^\ast$-representation of $L^1(G)$ on $\Hilbert$, i.e.\ a $^\ast$-homomorphism, $\tilde{\pi} \!: L^1(G) \to {\cal B}(\Hilbert)$ by letting
\[
  \langle \tilde{\pi}(f) \xi, \eta \rangle := \int_G f(x) \langle \pi(x) \xi, \eta \rangle \, dx \qquad (f \in L^1(G), \, \xi, \eta \in \Hilbert).
\]
Moreover, any $^\ast$-representation of $L^1(G)$ arises from a representation of $G$ in the above fashion (see \cite[\S 13]{Dix}).
\begin{examples}
\item The {\it left regular representation\/} $\lambda$ of $G$ on $L^2(G)$ is defined by 
\[
  (\lambda(x) \xi)(y) := \xi(x^{-1} y) \qquad (x,y \in G, \, \xi \in L^2(G)).
\]
\item Similarly, the {\it right regular representation\/} $\rho$ of $G$ on $L^2(G)$ is defined by 
\[
  (\rho(x) \xi)(y) := \frac{1}{\Delta(x)^\frac{1}{2}} \xi(yx) \qquad (x,y \in G, \, \xi \in L^2(G)).
\]
(Division by the square root of $\Delta(x)$ is necessary in order for
$\rho(x)$ to be an isometry.)
\item A function $\phi \!: G \to \comps$ is called positive definite if
\[
  \sum_{j,k=1}^n \overline{c_j} c_k \phi(x^{-1}_j x_k) \geq 0 \qquad (n \in \posints, \, c_1, \ldots, c_n \in \comps, \, x_1, \ldots, x_n \in G).
\]
Let $P(G)$ denote the continuous positive functions $\phi$ on $G$ such that $\phi(e) = 1$. The Gelfand--Naimark--Segal construction then yields a $^\ast$-representation of $L^1(G)$ --- and thus a representation $\pi_\phi$ of $G$ ---
on some Hilbert space $\Hilbert_\phi$.
\item Let $\Hilbert_u := \text{$\ell^2$-}\bigoplus_{\phi \in P(G)} \Hilbert_\phi$. Then $\pi_u := \bigoplus_{\phi \in P(G)} \pi_\phi$ is a representation of $G$ on $\Hilbert_u$, the {\it universal representation\/} of $G$. 
\end{examples}
\par
The {\it (full) group $\cstar$-algebra\/} of $G$ is the norm closure of $\tilde{\pi}_u(L^1(G))$ in ${\cal B}(\Hilbert_u)$ and denoted by $\cstar(G)$. Being a $\cstar$-algebra, $\cstar(G)$ is often easier to handle than $L^1(G)$, but this convenience 
comes at a price:
\begin{example}
The group $\cstar$-algebras of the groups $\ints/ 2\ints \times \ints/ 2\ints$ and $\ints/ 4\ints$ are $4$-dimensional, commutative $\cstar$-algebras and thus isometrically isomorphic to $\comps^4$. Nevertheless, 
$\ints/ 2\ints \times \ints/ 2\ints$ and $\ints/ 4\ints$ fail to be isomorphic.
\end{example}
\par
This stands in strong contrast to Theorem \ref{wendel1}. We will not deal so much with $\cstar(G)$ itself in this survey, but rather with its dual space (Sections \ref{FFSt} and \ref{QBan2}).
\par
Concluding this section, we briefly touch upon another Banach algebra associated with a locally compact group. 
\par
Given a locally compact group $G$, we denote by $M(G)$ the Banach space of all (finite) complex Borel measures on $G$. Via Riesz' 
representation theorem, $M(G)$ can be identified with the dual space of ${\cal C}_0(G)$, the space of all continuous functions on $G$ that vanish at infinity. We equip $M(G)$ with a convolution product through
\begin{equation} \label{conv2}
  \langle f, \mu \ast \nu \rangle := \int_{G \times G} f(xy) \, d\mu(x) \, d\nu(y) \qquad (\mu, \nu \in M(G), \, f \in {\cal C}_0(G)).
\end{equation}
This turns $M(G)$ into a Banach algebra (necessarily with identity).
\par
There are various important subspaces of $M(G)$. A measure $\mu \in M(G)$ is called {\it continuous\/} if $\mu( \{ x \} ) = 0$ for each $x \in G$. Let $M_c(G)$ denote the continuous measures in $M(G)$ and let $M_d(G)$ denote the
discrete measures. Then we have a direct sum decomposition
\begin{equation} \label{dirsum}
  M(G) = M_d(G) \oplus_{\ell^1} M_c(G).
\end{equation}
Of course, if $G$ is discrete, $M_c(G) = \{ 0 \}$ holds, so that $M(G) = M_d(G) = \ell^1(G) = L^1(G)$. Even though Haar measure need not be $\sigma$-finite there are versions of the Radon--Nikod\'ym theorem for regular Borel measures that
can be applied to measures absolutely continuous with respect to Haar measure. Let $M_a(G)$ denote the space of all such measures in $M(G)$; then an appropriate Radon--Nikod\'ym theorem allows us to identify $L^1(G)$ with $M_a(G)$. For
non-discrete $G$, we thus obtain a refinement of (\ref{dirsum}), namely
\[
  M(G) = M_d(G) \oplus_{\ell^1} M_s(G) \oplus_{\ell^1} M_a(G),
\] 
where $M_s(G)$ consists of those measures in $M_c(G)$ which are singular with respect to Haar measure. The subspaces $M_c(G)$ and $M_a(G)$ of $M(G)$ are in fact ideals, whereas $M_d(G)$ is only a subalgebra. Moreover, the identification
$L^1(G) \cong M_a(G)$ is an isometric isomorphism of Banach algebras, i.e.\ for $f, g \in L^1(G) = M_a(G)$, the two convolution formulae (\ref{conv}) and (\ref{conv2}) yield the same result.
\par
All of the above --- except the definition of $\cstar(G)$, which is covered in \cite{Dix} --- can be found in the encyclopedic treatise \cite{HR}. Less voluminous introductions to abstract harmonic analysis are \cite{Foll} and \cite{Rei}, which
has recently had an updated second edition (\cite{Rei2}). A monograph that solely focusses on the abelian case with its many peculiar features is \cite{Rud}.
\section{Amenable, locally compact groups}
Let $G$ be a locally compact group. A {\it mean\/} on $L^\infty(G)$ is a state of the commutative von Neumann algebra $L^\infty(G)$, i.e.\ a bounded linear functional $m \!: L^\infty(G) \to \comps$ such that $\| m \| = \langle 1, m \rangle = 1$. 
For $\phi \!: G \to \comps$, we define its left translate $L_x \phi$ by $x \in G$ through $(L_x \phi)(y) := \phi(xy)$ for $y \in G$. A mean on $L^\infty(G)$ is called {\it left invariant\/} if
\[
  \langle L_x \phi, m \rangle = \langle \phi, m \rangle \qquad (x \in G, \, \phi \in L^\infty(G)).
\]
\begin{definition} \label{amdef}
A locally compact group $G$ is called {\it amenable\/} if there is a left invariant mean on $L^\infty(G)$. 
\end{definition}
\par
Instead of via $L^\infty(G)$, the amenable locally compact groups can also be characterized through the existence of a left invariant mean on certain, much smaller subspaces of $L^\infty(G)$. In particular, $G$ is amenable if and only if
there is a left invariant mean on ${\cal C}_{\mathrm{b}}(G)$, the space of bounded continuous functions on $G$. Since ${\cal C}_{\mathrm{b}}(G) \subset \ell^\infty(G)$, this immediately yields that a locally compact group which is amenable
as a discrete group, is already amenable; the converse is false as we shall see below.
\par
The adjective ``a{\it men\/}able'' for the groups described in Definition \ref{amdef} was introduced by M.\ M.\ Day (\cite{Day}), apparently with a pun in mind: they have a me(a)n and at the same time are very tractable and thus truly amenable
in the sense of that adjective in colloquial English.
\par
If $G$ is compact, the inclusion $L^\infty(G) \subset L^1(G)$ holds trivially and Haar measure is a left invariant mean on $L^\infty(G)$; if $G$ is abelian, the Markov--Kakutani fixed point theorem can be used to obtain an invariant mean: 
consequently, all compact and all abelian groups are amenable. The easiest example of a non-amenable group is probably the free group in two generators:
\begin{example}
Let $\free_2$ denote the free group in two generators, say $a$ and $b$. Then each element of $\free_2$ is a reduced word over the alphabet $\{ a,b, a^{-1}, b^{-1} \}$. For any $x \in \{ a,b , a^{-1}, b^{-1} \}$ let
\[
  W(x) := \{ w \in \free_2: \text{$w$ starts with $x$} \}, 
\]
so that
\begin{equation} \label{free1}
  \free_2 = \{ e \} \cup W(a) \cup W(b) \cup W(a^{-1}) \cup W(b^{-1}),
\end{equation}
the union being disjoint. If $w \in \free_2 \setminus W(a)$, we necessarily have $a^{-1} w \in W(a^{-1})$ and thus $w \in a W(a^{-1})$; it follows that
\begin{equation} \label{free2}
  \free_2 = W(a) \cup a W(a^{-1}).
\end{equation}
Analogously, 
\begin{equation} \label{free3}
  \free_2 = W(b) \cup b W(b^{-1}).
\end{equation}
holds. Assume that there is a left invariant mean $m$ on $\ell^\infty(\free_2)$; it is not hard to see that $m$ must be positive, i.e.\ maps non-negative functions to non-negative numbers. Consequently, we obtain:
\begin{eqnarray*}
  1 & = & \langle 1, m \rangle \\
  & = & \langle \chi_{\{ e \}}, m \rangle + \langle \chi_{W(a)} , m \rangle + \langle \chi_{W(b)} , m \rangle + \langle \chi_{W(a^{-1})} , m \rangle + \langle \chi_{W(b^{-1})} , m \rangle, \quad\text{by (\ref{free1})}, \\
  & \geq & \langle \chi_{W(a)} , m \rangle + \langle \chi_{W(b)} , m \rangle + \langle \chi_{W(a^{-1})} , m \rangle + \langle \chi_{W(b^{-1})} , m \rangle, \quad\text{since $m$ is positive}, \\
  & \geq & \langle \chi_{W(a)} , m \rangle + \langle \chi_{W(b)} , m \rangle + \langle \chi_{aW(a^{-1})} , m \rangle + \langle \chi_{bW(b^{-1})} , m \rangle, \\
  &      & \quad\text{since $m$ is left invariant}, \\ 
  & = & \langle \chi_{W(a)} + \chi_{aW(a^{-1})}, m \rangle + \langle \chi_{W(b)} + \chi_{bW(b^{-1})} , m \rangle \\
  & \geq & \langle \chi_{W(a) \cup a W(a^{-1})}, m \rangle + \langle \chi_{W(b) \cup b W(b^{-1})}, m \rangle, \quad\text{again by the positivity of $m$}, \\
  & = & 1 + 1, \quad\text{by (\ref{free2}) and (\ref{free3})}, \\
  & = & 2.
\end{eqnarray*}
This, of course, is nonsense.
\end{example}
\par
One of the features making amenable groups genuinely amenable are their pleasant hereditary properties which we sum up in the following theorem:
\begin{theorem} \label{amher}
Let $G$ be a locally compact group.
\begin{items}
\item If $G$ is amenable and $H$ is a closed subgroup of $G$, then $H$ is amenable.
\item If $G$ is amenable, $H$ is another locally compact group, and $\theta \!: G \to H$ is a continuous homomorphism with dense range, then $H$ is amenable.
\item If $N$ is a closed, normal subgroup of $G$ such that both $N$ and $G/N$ are amenable, then $G$ is amenable.
\item If $( H_\alpha )_\alpha$ is an increasing family of closed subgroups of $G$ such that each $H_\alpha$ is amenable and such that $\bigcup_\alpha H_\alpha$ is dense in $G$, then $G$ is amenable.
\end{items}
\end{theorem} 
\par
Theorem \ref{amher} immediately increases our stock of both amenable and non-amenable, locally compact groups: All solvable as well as all locally finite groups are amenable (by Theorem \ref{amher}(iii) and (iv)) whereas every
locally compact group containing a copy of $\free_2$ as a closed subgroup cannot be amenable (by Theorem \ref{amher}(ii)); this yields, for example, the non-amenability of many Lie groups such as $\operatorname{SL}(\comps, N)$ for $N \geq 2$.
Also, Theorem \ref{amher}(ii) shows that an amenable, locally compact group need not be amenable as a discrete group: With a little linear algebra, it can be shown that the compact (and thus amenable) Lie group $\operatorname{SO}(3)$
contains an isomorphic copy of $\free_2$. Since in the discrete topology every subgroup is closed, $\operatorname{SO}(3)$ equipped with the discrete topology is not amenable.
\par
The hereditary properties listed in Theorem \ref{amdef} are all proven in \cite[Chapter 1]{LoA}.
\par
The standard reference for amenable, locally compact groups (and amenable semigroups, which we didn't define here) is the monograph \cite{Pat}. Older sources are \cite{Pie} and \cite{Gre}. 
\section{Homological properties of group algebras}
Let $\A$ a be a Banach algebra. Then a {\it left multiplier\/} of $\A$ is a bounded linear map $L \!: \A \to \A$ satisfying $L(ab) = aLb$ for all $a,b \in \A$. Often, it is possible to give concrete descriptions of the left multipliers
of a given Banach algebra as is done in the following result from \cite{Wen}:
\begin{proposition} \label{wendel2}
Let $G$ be a locally compact group, and let $L$ be a left multiplier of $L^1(G)$. Then there is $\mu \in M(G)$ such that
\[
  Lf = f \ast \mu \qquad (f \in L^1(G)).
\]
\end{proposition}
\begin{proof}
By Theorem \ref{BAI}, $L^1(G)$ has a bounded approximate identity, say $( e_\alpha )_\alpha$. The net $( L e_\alpha )_\alpha$ is bounded in the dual space $M(G) \cong {\cal C}_0(G)^\ast$ and thus has a $w^\ast$-accumulation point $\mu$.
This $\mu$ works.
\end{proof}
\par
Proposition \ref{wendel2} can be interpreted in terms of Hochschild cohomology.
\par
We will not attempt to define Hochschild cohomology groups of arbitrary order, but confine ourselves to first Hochschild cohomology groups. For more, see \cite{Joh1} and \cite[Chapters 2 and 5]{LoA}.
\par
A bimodule $E$ over a Banach algebra $\A$ is called a {\it Banach $\A$-bimodule\/} if it is also a Banach space such that the module actions are continuous. A {\it derivation\/} from $\A$ into $E$ is a bounded linear map
$D \!: \A \to E$ satisfying
\[
  D(ab) = a \cdot Db + (Da) \cdot b \qquad (a,b \in \A);
\] 
we write ${\cal Z}^1(\A,E)$ for the Banach space of all derivations from $\A$ into $E$. A derivation $D$ is called {\it inner\/} if there is $x \in E$ such that
\[
  Da = a \cdot x - x \cdot a \qquad (a \in \A);
\]
in this case, we say that $x$ {\it implements\/} $D$. The space of all inner derivations is denoted by ${\cal B}^1(\A,E)$.
\begin{definition}
Let $\A$ be a Banach algebra, and let $E$ be a Banach $\A$-bimodule. Then the {\it first Hochschild cohomology group of $\A$ with coefficients in $E$\/} is defined as
\[
  {\cal H}^1(\A,E):= {\cal Z}^1(\A,E) / {\cal B}^1(\A,E).
\]
\end{definition}
\par
Note that the quotient topology on ${\cal H}^1(\A,E)$ need not be Hausdorff.
\par
Hochschild cohomology is named in the honor of G.\ Hochschild who introduced it in the 1940s (\cite{Hoch1} and \cite{Hoch2}) --- in a purely algebraic context, of course. The first to adapt it to the Banach algebra context was H.\ Kamowitz
(\cite{Kam}).
\begin{example}
Let $G$ be a locally compact group, let $\A := L^1(G)$, and let $E := M(G)$ be equipped with the following module operations:
\[
  f \cdot \mu := f \ast \mu \quad\text{and}\quad \mu \cdot f := 0 \qquad (f \in L^1(G), \, \mu \in M(G)).
\]
Then ${\cal Z}^1(\A,E)$ consists precisely of the left multipliers of $L^1(G)$, and Proposition \ref{wendel2} is equivalent to the assertion that ${\cal H}^1(\A,E) = \{ 0 \}$.
\end{example}
\par
Of course, a much more natural module action of $L^1(G)$ on $M(G)$ is via convolution from the right and from the left. The following problem is therefore rather natural in a cohomological context and was a main reason for
B.\ E.\ Johnson to develop his theory of amenable Banach algebras:
\begin{problem}
Let $G$ be a locally compact group. Does ${\cal H}^1(L^1(G),M(G)) = \{ 0 \}$ hold or --- equivalently --- is there, for each derivation $D \!: L^1(G) \to L^1(G)$, a measure $\mu \in M(G)$ such that
\[
  Df = f \ast \mu - \mu \ast f \qquad (f \in L^1(G))
\] 
holds?
\end{problem}
\par
We would like to mention that this problem was solved affirmatively only recently in its full generality by V.\ Losert --- after having been open for more than three decades.
\par
To connect the problem with Hochschild cohomology (and give an affirmative answer for amenable $G$), we need one more definition. Given a Banach algebra $\A$ and an Banach $\A$-bimodule $E$, the dual space $E^\ast$ of $E$ becomes
a Banach $\A$-bimodule via
\[
  \langle a \cdot \phi, x \rangle := \langle \phi, x \cdot a \rangle \quad\text{and}\quad \langle \phi \cdot a , x \rangle := \langle \phi, a \cdot x \rangle
  \qquad (a \in \A, \, \phi \in E^\ast, \, x \in E).
\]
The following theorem due to B.\ E.\ Johnson (\cite[Theorem 2.5]{Joh1}) is the starting point of the theory of amenable Banach algebras:
\begin{theorem} \label{Barry}
Let $G$ be a locally compact group. Then the following are equivalent:
\begin{items}
\item $G$ is amenable.
\item ${\cal H}^1(L^1(G),E^\ast) = \{ 0 \}$ for each Banach $L^1(G)$-bimodule $E$.
\end{items}
\end{theorem}
\begin{proof}
To keep matters simple, we only treat the discrete case, so that, in particular, $L^1(G) = \ell^1(G)$ and $L^\infty(G) = \ell^\infty(G)$.
\par
(i) $\Longrightarrow$ (ii): Let $E$ be a Banach $\ell^1(G)$-bimodule, and let $D \!: \ell^1(G) \to E^\ast$ be a derivation. For $x \in E$, define $x_\infty \in \ell^\infty(G)$ through
\[
  x_\infty(g) := \langle x, \delta_g \cdot D(\delta_{g^{-1}}) \rangle \qquad (g \in G).
\]
Let $m$ be a left invariant mean on $\ell^\infty(G)$, and define $\phi \in E^\ast$, by letting
\[
  \langle x, \phi \rangle := \langle x_\infty, m \rangle \qquad (x \in E).
\]
It is routinely, albeit a bit tediously verified that
\[
  Df = f \cdot \phi - \phi \cdot f \qquad (f \in \ell^1(G)).
\]
\par
(ii) $\Longrightarrow$ (i): Turn $\ell^\infty(G)$ into a Banach $\ell^1(G)$-module by letting
\[
  f \cdot \phi := f \ast \phi \quad\text{and}\quad \phi \cdot f := \left( \sum_{x \in G} f(x) \right) \phi \qquad (f \in \ell^1(G), \, \phi \in \ell^\infty(G)).
\]
Note that $\comps 1$ is a submodule of $\ell^\infty(G)$, so that it makes sense to define $E := \ell^\infty(G) / \comps 1$. Fix $n \in \ell^\infty(G)^\ast$ with $\langle 1, n \rangle = 1$. It is easily seen that the inner derivation
$D \!: \ell^1(G) \to \ell^\infty(G)^\ast$ implemented by $n$ maps, in fact, into $E^\ast$. Hence, by (ii), there is $\tilde{n} \in E^\ast$ implementing $D$. Letting
\[
  m := \frac{|n - \tilde{n}|}{\| n - \tilde{n} \|}
\]
we obtain a left invariant mean on $\ell^\infty(G)$.
\end{proof}
\begin{corollary}
Let $G$ be a locally compact group. Then ${\cal H}^1(L^1(G), M(G)) = \{ 0 \}$ holds whenever $G$ is amenable.
\end{corollary}
\par
The choice of adjective in the following definition should be clear in view of Theorem \ref{Barry}:
\begin{definition} \label{bamdef}
A Banach algebra $\A$ is said to be {\it amenable\/} if ${\cal H}^1(\A,E^\ast) = \{ 0 \}$ for each Banach $\A$-bimodule $E$.
\end{definition}
\par
The following is an elementary, but useful property of amenable Banach algebras:
\begin{proposition} \label{bambai}
Let $\A$ be an amenable Banach algebra. Then $\A$ has a bounded approximate identity.
\end{proposition}
\begin{proof}
Let $\A$ be equipped with the module actions
\[
  a \cdot x := ax \quad\text{and}\quad x \cdot a := 0 \qquad (a, x \in \A).
\]
Then the canonical inclusion of $\A$ in $\A^{\ast\ast}$ is a derivation and thus inner, i.e.\ there is $E \in \A^{\ast\ast}$ such that $a \cdot E = a$ for $a \in \A$. Let $( e_\alpha )_\alpha$ be a bounded net in $\A$ that
converges to $E$ in the $w^\ast$-topology of $\A^{\ast\ast}$; it follows $a e_\alpha \to a$ in the weak topology for all $a \in \A$. Passing to convex combinations, we can achieve that $a a_\alpha \to a$ for all $a \in \A$
in the norm topology. Hence, $\A$ has a bounded left approximate identity.
\par
Analogously, one shows that $\A$ has also a bounded right approximate identity. But then $\A$ already has a (two-sided) bounded approximate identity (\cite[Proposition 2.9.3]{Dal})
\end{proof}
\par
The theory of amenable Banach algebra has been an active, ever expanding area of research since its inception in \cite{Joh1}. For $\cstar$-algebras, amenability in the sense of Definition \ref{bamdef} is equivalent to the important
property of nuclearity; see \cite[Chapter 6]{LoA} for a self-contained exposition of this equivalence.
\par
A drawback of Definition \ref{bamdef} is that it is based on a condition for {\it all\/} Banach bimodules over a given algebra. For practical purposes, i.e.\ to confirm or to rule out whether or not a given Banach algebra is amenable,
it is therefore often difficult to handle. There is, however, a more intrinsic characterization of amenable Banach algebra, which is also due to Johnson (\cite{Joh2}).
\par
Following \cite{ER}, we denote the (completed) projective tensor product of two Banach spaces by $\tensor^\gamma$. Given a Banach algebra $\A$ the tensor product $\A \tensor^\gamma \A$ becomes a Banach $\A$-bimodule via
\[
  a \cdot (x \tensor y) := ax \tensor y \quad\text{and}\quad (x \tensor y) \cdot a := x \tensor ya \qquad (a,x, y \in\A).
\]
Multiplication induces a bounded linear map $\Gamma \!: \A \tensor^\gamma \A \to \A$ which is easily seen to be an $\A$-bimodule homomorphism.
\par
The following equivalence is from \cite{Joh2}:
\begin{proposition} \label{diagprop}
The following are equivalent for a Banach algebra $\A$:
\begin{items}
\item $\A$ is amenable.
\item There is an approximate diagonal for $\A$, i.e.\ a bounded net $(\mathbf{m}_\alpha)_\alpha$ in $\A \tensor^\gamma \A$ such that
\[
  a \cdot \mathbf{m}_\alpha - \mathbf{m}_\alpha \cdot a \to 0 \quad\text{and}\quad a \Gamma \mathbf{m}_\alpha \to a \qquad (a \in \A).
\]
\item There is a virtual diagonal for $\A$, i.e.\ an element $\mathbf{M} \in (\A \tensor^\gamma \A)^{\ast\ast}$ such that
\[
  a \cdot \mathbf{M} = \mathbf{M} \cdot a \quad\text{and}\quad a \Gamma^{\ast\ast} \mathbf{M} = a \qquad (a \in \A).
\]
\end{items}
\end{proposition} 
\begin{proof}
Every $w^\ast$-accumulation point of an approximate diagonal is a virtual diagonal: this settles (ii) $\Longrightarrow$ (iii). The converse is a simple approximation argument.
\par
For the equivalence of (i) and (ii), we suppose for the sake of simplicity that $\A$ has an identity $e$.
\par
(i) $\Longrightarrow$ (iii): Let
\[
  D \!: \A \to \ker \Gamma^{\ast\ast}, \quad a \mapsto a \tensor e - e \tensor a.
\]
Then $D$ is inner, i.e.\ there is $\mathbf{N} \in \ker \Gamma^{\ast\ast}$ implementing it. Hence, $\mathbf{M} := e \tensor e - \mathbf{N}$ is a virtual diagonal.
\par
(ii) $\Longrightarrow$ (i): Let $E$ be a Banach $\A$-bimodule, and let $D \!: \A \to E^\ast$ be a derivation. We can confine ourselves to the case when $E$, and thus $E^\ast$, is unital. Let $( \mathbf{m}_\alpha )_\alpha$ be an approximate
diagonal for $\A$, and let $\phi \in E^\ast$ a $w^\ast$-accumulation point of $(\Gamma_{E^\ast}((\id_\A \tensor D)\mathbf{m}_\alpha))_\alpha$, where $\Gamma_{E^\ast} \!: \A \tensor^\gamma E^\ast \to E^\ast$ is the linear
map induced by the left module action of $\A$ on $E^\ast$. Then $\phi$ implements $D$.
\end{proof}
\par
Besides being more concrete then Definition \ref{bamdef}, Proposition \ref{diagprop}(ii) and (iii) allow a refinement of the notion of amenability: $\A$ is $C$-amenable with $C \geq 1$ if is has a virtual diagonal of norm at most $C$.
\par
In analogy with Theorem \ref{amher}, amenability for Banach algebras has nice hereditary properties:
\begin{theorem} \label{bamher}
Let $\A$ be a Banach algebra.
\begin{items}
\item If $\A$ is amenable, $\B$ is another Banach algebra, and $\theta \!: \A \to \B$ is a continuous homomorphism with dense range, then $\B$ is amenable; in particular, $\A/I$ is amenable for every
closed ideal $I$ of $\A$.
\item If $I$ is a closed ideal of $\A$ such that both $I$ and $\A / I$ are amenable, then $\A$ is amenable.
\item If $I$ is a closed ideal of $\A$, then $I$ is amenable if and only if it has a bounded approximate identity and if and only if its annihilator $I^\perp$ in $\A^\ast$ is complemented in $\A^\ast$.
\item If $( \A_\alpha )_\alpha$ is a directed family of closed subalgebras of $\A$ such that each $\A_\alpha$ is $C$-amenable for some universal $C \geq 1$ and such that $\bigcup_\alpha \A_\alpha$ is dense in $\A$, then $\A$ is amenable.
\end{items}
\end{theorem}
\par
For proofs, see \cite[Chapter 2]{LoA}, for instance.
\par
By replacing the class of dual bimodules in Definition \ref{bamdef}, one can, of course, weaken or strengthen the notion of amenability. We limit ourselves to looking at only one of those variants:
\begin{definition} \label{wamdef}
A Banach algebra $\A$ is called {\it weakly amenable\/} if ${\cal H}^1(\A,\A^\ast) = \{ 0 \}$.
\end{definition}
\par
Weak amenability was introduced by W.\ G.\ Bade, P.\ C.\ Curtis, Jr., and H.\ G.\ Dales in \cite{BCD} for commutative Banach algebra (using a formally stronger, but in fact equivalent condition). Definition \ref{wamdef}, as we use it, 
originates in \cite{JohWA1}.
\par
To illustrate how much weaker than amenability weak amenability is, we state the following theorem due to Johnson (\cite{JohWA2}) and sketch its ingeniously simple proof by M.\ Despi\'c and F.\ Ghahramani (\cite{DG}).
\begin{theorem}
Let $G$ be a locally compact group. Then $L^1(G)$ is weakly amenable.
\end{theorem}
\begin{proof}
We only consider the discrete case. We may identify $\ell^1(G)^\ast$ with $\ell^\infty(G)$; let $D \!: \ell^1(G) \to \ell^\infty(G)$ be a derivation. Since $\ell^\infty_\reals(G)$, the space of all $\reals$-valued bounded functions on $G$, is a complete 
lattice, the function
\[
  \phi := \sup \{ \re \, (D\delta_x) \cdot \delta_{x^{-1}} : x \in G \} + i \,  \sup \{ \im \, (D\delta_x) \cdot \delta_{x^{-1}} : x \in G \}  
\]
exists and lies in $\ell^\infty(G)$. It is easily seen to implement $D$.
\end{proof}
\par
So far, we have only considered $L^1(G)$ in this section. We now turn to $M(G)$.
\par
The following recent theorem due to Dales, Ghahramani, and Helemski\u{\i} (\cite{DGH}), characterizes those locally compact groups $G$, for which $M(G)$ is weakly amenable and amenable, respectively:
\begin{theorem} \label{DGHthm}
Let $G$ be a locally compact group. Then $M(G)$ is weakly amenable if and only if $G$ is discrete. In particular, $M(G)$ is amenable if and only if $G$ is discrete and amenable.
\end{theorem}
\par
The fairly intricate proof centers around showing that the closed linear span of $\{ \mu \ast \nu : \mu, \nu \in M_c(G)\}$ has infinite codimension in $M_c(G)$. The technical heart of the argument is the construction --- in the
metrizable case --- of a perfect subset $V$ of $G$ that supports a continuous measure, but such that $(\mu \ast \nu)(V) = 0$ for all $\mu,\nu \in M_c(G)$.
\par
Homological algebra can be systematically equipped with functional analytic overtones: Attempts in this direction were made by several mathematicians in the 1960s and early 1970s --- most persistently by Helemski\u{\i} and his
Moscow school (see the monograph \cite{Hel} or the more introductory text \cite{Hel2} or the even more introductory \cite[Chapter 5]{Hel}). A central r\^ole in Helemski\u{i}'s approach is played by the notion of a projective module:
\begin{definition} 
Let $\A$ be a Banach algebra. A Banach $\A$-bimodule $E$ is called {\it projective\/} if, for each Banach $\A$-bimodule $F$ and each bounded $\A$-bimodule homomorphism $\pi \!: E \to F$ with a bounded linear right inverse, there
is a bounded $\A$-bimodule homomorphism $\rho \!: F \to E$ such that $\pi \circ \rho = \id_F$.
\end{definition}
\par
We suppose the existence of merely a {\it linear\/} right inverse, and projectivity gives us a right inverse that respects the module actions.
\begin{definition} 
A Banach $\A$ is called {\it biprojective\/} if it is a projective Banach $\A$-bimodule.
\end{definition}
\par
We quote the following characterization of biprojective Banach algebras without proof:
\begin{proposition} \label{biproj}
A Banach algebra $\A$ is biprojective if and only if $\Gamma \!: \A \tensor^\gamma \A \to \A$ has a bounded right inverse which is an $\A$-bimodule homomorphism.
\end{proposition}
\par
The following is \cite[Theorem 51]{HelTMMS}:
\begin{theorem} \label{sasha0}
The following are equivalent for a locally compact group:
\begin{items}
\item $G$ is compact.
\item $L^1(G)$ is biprojective.
\end{items}
\end{theorem}
\begin{proof}
(i) $\Longrightarrow$ (ii): We identify $L^1(G) \tensor^\gamma L^1(G)$ and $L^1(G \times G)$. Define $\rho \!: L^1(G) \to L^1(G \times G)$ by letting
\[
  \rho(f)(x,y) := f(xy) \qquad (f \in L^1(G), \, x,y \in G).
\]
Then $\rho$ is a right inverse of $\Gamma$ as required by Proposition \ref{biproj}.
\par
(ii) $\Longrightarrow$ (i): The {\it augmentation character\/}
\begin{equation} \label{augment}
  L^1(G) \to \comps, \quad f \mapsto \int_G f(x) \, dx
\end{equation}
turns $\comps$ into a Banach $L^1(G)$-bimodule. Since $\comps$ is a quotient of $L^1(G)$ one can show that $\comps$ must also be a projective Banach $L^1(G)$-bimodule. Hence, there is a bounded right inverse $\rho$ of (\ref{augment}) which is
also an $L^1(G)$-bimodule homomorphism. It is easy to see that $\rho(1) \in L^1(G)$ must be translation invariant and therefore constant. This is possible only if $G$ is compact.
\end{proof}
\par
A notion equally central to topological homology as projectivity is
that of flatness. A Banach algebra $\A$ is (obviously) called {\it
  biflat\/} if it is a flat Banach bimodule over itself. We shall not define here flat Banach bimodules in general,
but use an equivalent condition (similar to Proposition \ref{biproj}) in order to introduce biflat Banach algebras:
\begin{definition} \label{biflat}
A Banach algebra $\A$ is called {\it biflat\/} if there is a bounded $\A$-bimodule homomorphism $\rho \!: \A \to (\A \tensor^\gamma \A)^{\ast\ast}$ such that $\Gamma^{\ast\ast} \circ \rho$ is the canonical embedding of $\A$ into $\A^{\ast\ast}$.
\end{definition}
\par
Obviously, biflatness is weaker than biprojectivity.
\par
The following result (also due to Helemski\u{\i}) relates biflatness and amenability:
\begin{theorem} \label{sasha}
The following are equivalent for a Banach algebra $\A$:
\begin{items}
\item $\A$ is amenable.
\item $\A$ is biflat and has a bounded approximate identity.
\end{items}
\end{theorem}
\begin{proof}
(i) $\Longrightarrow$ (ii): By Proposition \ref{bambai}, amenable Banach algebras always have a bounded approximate identity. Let $\mathbf{M} \in (\A \tensor^\gamma \A)^{\ast\ast}$ be a virtual diagonal for $\A$. Then
\[
  \rho \!: \A \to  (\A \tensor^\gamma \A)^{\ast\ast}, \quad a \mapsto a \cdot \mathbf{M}
\]
is a bimodule homomorphism as required in Definition \ref{biflat}.
\par
(ii) $\Longrightarrow$ (i): Let $\rho  \!: \A \to  (\A \tensor^\gamma \A)^{\ast\ast}$ an $\A$-bimodule homomorphism as in Definition \ref{biflat} and let $( e_\alpha )_\alpha$ be a bounded approximate identity for $\A$.
Then any $w^\ast$-accumulation point of $( \rho(e_\alpha) )_\alpha$ is a virtual diagonal for $\A$.
\end{proof}
\par
Combining this with Theorem \ref{BAI}, we see that $L^1(G)$ is biflat if and only if $L^1(G)$ is amenable, i.e.\ if and only if $G$ is amenable.
\par
As mentioned several times already, he theory of amenable Banach algebras was initiated in \cite{Joh1}. Recent expositions can be found in \cite{Dal} and \cite{LoA}. For Helemski\u{\i}'s approach to topological homology, see his books \cite{Hel} and 
\cite{Hel2}, and also the survey article \cite{HelTMMS}.
\section{Fourier and Fourier--Stieltjes algebras} \label{FFSt}
The {\it dual group\/} or {\it character group\/} of an abelian, locally compact group $G$ is the set of all continuous group homomorphisms from $G$ into $\torus$; equipped with pointwise multiplication and the compact open topology, 
it becomes a locally compact group in its own right, which we denote by $\Hat{G}$. The {\it Fourier--Stieltjes transform\/} ${\cal FS} \!: M(G) \to {\cal C}_{\mathrm{b}}(G)$ is defined via
\[
  {\cal FS}(\mu)(\gamma) := \int_G \overline{\gamma(x)} \, d\mu(x) \qquad (\mu \in M(G), \, \gamma \in \Hat{G}).
\]
It is a continuous, injective homomorphism of Banach algebras. The range of ${\cal FS}$ is denoted by $B(\Hat{G})$ and called the {\it Fourier--Stieltjes algebra\/} of $\Hat{G}$; by definition, it is isometrically isomorphic to $M(G)$.
\par
The restriction $\cal F$ of $\cal FS$ to $L^1(G)$ is called the {\it Fourier transform\/}. The Riemann--Lebesgue lemma yields immediately that the range of $\cal F$ is contained in ${\cal C}_0(\hat{G})$. We call ${\cal F}(L^1(G))$ the {\it
Fourier algebra\/} of $\Hat{G}$ and denote it by $A(\Hat{G})$.
\par
Since $\Hat{\Hat{G}} \cong G$ for every abelian, locally compact group $G$, the algebras $A(G)$ and $B(G)$ are defined for {\it every\/} such group.
\par
To extend the definitions of $A(G)$ and $B(G)$ to arbitrary --- not necessarily abelian --- locally compact groups, first recall that (\cite[Theorem 1.6.3]{Rud})
\begin{equation} \label{Fourier}
  A(G) = \{ \xi \ast \eta : \xi, \eta \in L^2(G) \}
\end{equation}
for each abelian, locally compact group $G$: this follows from Plancherel's theorem and the elementary fact that each function in $L^1(\Hat{G})$ is the pointwise product of two $L^2$-functions (the convolution on the right hand side of
(\ref{Fourier}) is formally defined as in (\ref{conv})). When trying to use the right hand side of (\ref{Fourier}) to define $A(G)$ for arbitrary $G$, we are faced with the problem that the convolution product of two $L^2$-functions need not
even be defined (unless the modular function of $G$ is trivial as in the abelian case). We define for an arbitrary function $f \!: G \to \comps$ another function $\check{f} \!: G \to \comps$ defined by $\check{f}(x) := f(x^{-1})$ for $x \in G$.
If $\Delta \equiv 1$, the map $L^2(G) \ni \xi \mapsto \check{\xi}$ is a unitary operator of $L^2(G)$; otherwise, it may not even leave $L^2(G)$ invariant. Anyway, the convolution product $\xi \ast \check{\eta}$ is well-defined and lies in
${\cal C}_0(G)$ for all $\xi, \eta \in L^2(G)$.
\par
We can therefore {\it define\/} the {\it Fourier algebra\/} of $G$ by letting
\begin{equation} \label{Fourier2}
  A(G) := \{ \xi \ast \check{\eta} : \xi, \eta \in L^2(G) \};
\end{equation}
we equip it with a norm via
\[
  \| f \|_{A(G)} := \inf \{ \| \xi \| \| \eta \| : \xi, \eta \in L^2(G), \, f = \xi \ast \check{\eta} \} \qquad (f \in A(G)).
\]
\par
Admittedly, it is not evident from (\ref{Fourier2}) that $A(G)$ is an algebra (or even a linear space). Nevertheless, the following is true (compiled from \cite{Eym}):
\begin{proposition}
Let $G$ be a locally compact group. Then $A(G)$ is a regular, Tauberian, commutative Banach algebra whose character space is canonically identified with $G$.
\end{proposition}
\par
The Fourier algebra can be conveniently described in terms of the left regular representation $\lambda$ of $G$ (introduced in Section \ref{aha}) --- this will become particularly relevant in Section \ref{QBan} below: A function
$f \!: G \to \comps$ belongs to $A(G)$ if and only if there are $\xi, \eta \in L^2(G)$ such that 
\begin{equation} \label{Fourier3}
  f(x) = \langle \lambda(x) \xi, \eta \rangle \qquad (x \in G);
\end{equation}
we call such functions {\it coefficient functions\/} of $\lambda$. More generally, we call a coefficient function $f$ of a representation $\pi$ of $G$ on some Hilbert space $\Hilbert$ if there are $\xi,\eta \in \Hilbert$ such that
\begin{equation} \label{coeff}
  f(x) = \langle \pi(x) \xi, \eta \rangle \qquad (x \in G).
\end{equation}
We then define the {\it Fourier--Stieltjes algebra\/} of $G$ as
\begin{equation} \label{FSTdef}
  B(G) := \{ f : \text{$f$ is a coefficient function of a representation of $G$} \}.
\end{equation}
It is somewhat easier than for $A(G)$ to see that $B(G)$ is indeed an algebra: sum and products of functions correspond to direct sums and tensor products of representation; it can be equipped with a norm by letting
\[
  \| f \|_{B(G)} := \inf \{ \| \xi \| \| \eta \| : \text{$f$ is represented as in (\ref{coeff})} \} \qquad (f \in B(G)).
\]
\par
The following is again a summary of results from \cite{Eym}:
\begin{proposition}
Let $G$ be a locally compact group. Then $B(G)$ is a commutative Banach algebra with identity which contains $A(G)$ as a closed ideal.
\end{proposition}
\par
To see that $B(G)$ as defined in (\ref{FSTdef}) is the same as ${\cal FS}(M(\Hat{G}))$, note that every coefficient of the form $f(x) = \langle \pi(x) \xi, \xi \rangle$ for $x \in G$ with $\pi$ being a representation of $G$ on some
Hilbert space $\Hilbert$ containing $\xi$ is positive definite. Moreover, every continuous, positive definite function on $G$ arised in this fashion (GNS-construction) Hence, $B(G)$ is the linear span of the continuous, positive definite functions 
on $G$ but the same is true for ${\cal FS}(M(\Hat{G}))$ if $G$ is abelian (\cite[Bochner's theorem, 1.4.3]{Rud}).
\par
The first to characterize properties of $G$ in terms of $A(G)$ was H.\ Leptin in \cite{Lep}. He proved the following theorem, whose proof we omit:
\begin{theorem} \label{lep}
The following are equivalent for a locally compact group $G$:
\begin{items}
\item $G$ is amenable.
\item $A(G)$ has a bounded approximate identity.
\end{items}
\end{theorem}
\par
With Proposition \ref{bambai} in mind, we see at once that the amenability of $A(G)$ forces $G$ to be an amenable locally compact group. The tempting conjecture that the converse is true as well, however, is wrong
(\cite{JohFourier}):
\begin{theorem} \label{Barry2}
Let $G$ be an infinite, compact group which, for each $n \in \posints$, has only finitely many irreducible unitary representations. Then $A(G)$ is not amenable.
\end{theorem}
\par
Examples for such groups are, for instance, $\operatorname{SO}(N)$ for $N \geq 3$. For $G = \operatorname{SO}(3)$, the Fourier algebra $A(G)$ is not even weakly amenable (see also \cite{JohFourier}).
\par
We will not even outline a proof for Theorem \ref{Barry2} because we'll obtain a much stronger result in Section \ref{QBan2} below.
\par
On the positive side, we have the following result from \cite{LLW}:
\begin{theorem} \label{LLWthm}
Let $G$ be a locally compact group which has an abelian subgroup of finite index. Then $A(G)$ is amenable.
\end{theorem}
\begin{proof}
Let $H$ be an abelian subgroup of $G$ such that $[G :H ] =: N <\infty$. Without loss of generality, suppose that $H$ is closed (otherwise, replace it by its closure). 
It follows that $A(G) \cong A(H)^N$. Since $A(H) \cong L^1(\Hat{H})$ is amenable by Theorem \ref{Barry}, the hereditary properties of amenability yield the amenability of $A(G)$.
\end{proof}
\par
With operator space methods, we shall see in Section \ref{QBan2} below that the rather restrictive sufficient condition of Theorem \ref{LLWthm} to ensure the amenability of $A(G)$ is, in fact, necessary.
\par
All the facts about Fourier and Fourier--Stieltjes algebras of locally compact, abelian groups mentioned in this section are contained in \cite[Chapter 1]{Rud}; see also \cite{HR}, \cite{Rei}, \cite{Rei2}, and \cite{Foll}. For $A(G)$ and $B(G)$
with $G$ arbitrary, P.\ Eymard's seminal paper \cite{Eym} still seems to be the best reference.
\section{Operator spaces}
Naively one may think that the term ``operator space'' just designates (certain) spaces of bounded linear operators. Indeed, this used to be the case:
\begin{definition} \label{opdef1}
A {\it concrete operator space\/} is a closed subspace of ${\cal B}(\Hilbert)$ for some Hilbert space $\Hilbert$.
\end{definition}
\par
Given a Banach space $E$, one can easily construct an isometry from $E$ into a commutative $\cstar$-algebra ${\cal C}(\Omega)$ of continuous functions on some compact Hausdorff space
$\Omega$. Since ${\cal C}(\Omega)$ can be represented on some Hilbert space $\Hilbert$, the Banach space $E$ is isometrically isomorphic to some concrete operator space: In Definition \ref{opdef1}, it is not important
that a given space can somehow be found sitting in ${\cal B}(\Hilbert)$, but {\it how\/} it sits there. Operator spaces are sometimes refereed to as ``quantized Banach spaces''. This has not so much to do with their potential applications to
quantum physics but with a formal analogy: the observables in classical physics are functions whereas those in quantum physics are operators on Hilbert space. The process of replacing functions by operators is therefore often referred to
as ``quantization''.  Banach spaces, i.e.\ spaces of functions, thus belong into the ``classical'' realm whereas operator spaces are their quantized counterpart.
\par
The adjective ``concrete'' Definition \ref{opdef1} suggests that there may also be ``abstract'' operator spaces. To define them, we first have to introduce some notation.
\par
Given a linear space $E$ and $n,m \in \posints$, we write $M_{n,m}(E)$ to denote the space of $n \times m$ matrices with entries from $E$; if $n = m$, we simply write $M_n(E)$. For the sake of simplicity, we only write $M_{n,m}$
instead of $M_{n,m}(\comps)$ or even $M_n$ if $n=m$. Identifying $M_{n,m}$ with ${\cal B}(\ell^2_n,\ell^2_m)$, we equip $M_{n,m}$ with a norm which we denote by $\| \cdot \|$ throughout.
\begin{definition} \label{opdef2}
An {\it operator space\/} is a linear space $E$ with a complete norm $\| \cdot \|_n$ on $M_n(E)$ for each $n \in \posints$ such that
\begin{equation} \tag{R 1} \label{R1}
  \left\| \begin{array}{c|c} x & 0 \\ \hline 0 & y \end{array} \right\|_{n+m} = \max \{ \| x \|_n, \| y \|_m \}
  \qquad (n,m \in \posints, \, x \in M_n(E), \, y \in M_m(E))
\end{equation}
and
\begin{equation} \tag{R 2} \label{R2}
  \| \alpha  x  \beta \|_n \leq \| \alpha \| \| x \|_n \| \beta \| \qquad (n \in \posints, \, x \in M_n(E), \, \alpha, \beta \in M_n).
\end{equation}
\end{definition}
\begin{example}
Let $\Hilbert$ and $\mathfrak K$ be Hilbert spaces. For each $n \in \posints$, identify $M_n({\cal B}(\Hilbert,{\mathfrak K}))$ with ${\cal B}(\ell^2_n(\Hilbert), \ell^2_n({\mathfrak K}))$. The operator norm on each matrix level
${\cal B}(\ell^2_n(\Hilbert), \ell^2_n({\mathfrak K}))$ then turns ${\cal B}(\Hilbert,{\mathfrak K})$ (and each of its closed subspaces) into an operator space. In particular, each concrete operator space is an operator space 
in the sense of Definition \ref{opdef2}.
\end{example}
\par
Every concrete operator space is an operator space, but what about the converse? It is clear that this question can only be answered up to (the appropriate notion of) isomorphism.
\par
Given two linear spaces $E$ and $F$, a linear map $T \!: E \to F$, and $n \in \posints$, we define the the {\it $n$-th amplification\/} $T^{(n)} \!: M_n(E) \to M_n(F)$ by applying $T$ to each matrix entry.
\begin{definition} \label{cbdef}
Let $E$ and $F$ be operator spaces, and let $T \in {\mathcal B}(E,F)$. Then:
\begin{alphitems}
\item $T$ is {\it completely bounded\/} if
\[
  \| T \|_{\mathrm{cb}} := \sup_{n \in \posints} \left\| T^{(n)} \right\|_{{\cal B}(M_n(E),M_n(F))} < \infty.
\]
\item $T$ is a {\it complete contraction\/} if $\| T \|_{\mathrm{cb}} \leq 1$.
\item $T$ is a {\it complete isometry\/} if $T^{(n)}$ is an isometry for each $n \in \posints$.
\end{alphitems}
We denote the completely bounded operators from $E$ to $F$ by $\CB(E,F)$.
\end{definition}
\par
It is easy to see that $\CB(E,F)$ equipped with $\| \cdot \|_\cb$ is a Banach space.
\begin{example}
Let $\A$ and $\B$ be $\cstar$-algebras, and let $\pi \!: \A \to \B$ be a $^\ast$-homomorphism. Since $\pi^{(n)}$ is a $^\ast$-homomorphism as well for each $n \in \posints$ and therefore contractive, it follows that
$\pi$ is a complete contraction. 
\end{example}
\par
The following theorem due to Z.-J.\ Ruan (\cite{Rua}) marks the beginning of abstract operator space theory:
\begin{theorem} \label{Ruan1}
Let $E$ be an operator space. Then $E$ is isometrically isomorphic to a concrete operator space.
\end{theorem}
\par
To appreciate Theorem \ref{Ruan1}, one should think of it as the operator space analog of the aforementioned fact that every Banach space is isometrically isomorphic to a closed subspace of some commutative $\cstar$-algebra: One could 
use it to {\it define\/} Banach spaces. With this ``definition'', however, even checking, e.g., that $\ell^1$ is a Banach space or that quotients and dual spaces of Banach spaces are again Banach spaces is difficult if not impossible.
\par
The following are examples of operator spaces:
\begin{enumerate}
\item Let $E$ be any Banach space. Then $E$ can be embedded into a commutative $\cstar$-algebra. This defines an operator space structure over $E$. This operator space is called the {\it minimal operator space\/} over $E$ and
is denoted by $\MIN(E)$; it is independent of the concrete embedding --- all that matters is that the $\cstar$-algebra is commutative. The adjective ``minimal'' is due to the following fact: Given another operator space $F$, every
operator in ${\cal B}(F,E)$ lies in $\CB(F,\MIN(E))$ such that its $\cb$-norm is just the operator norm. 
\item Given any Banach space $E$, we define an operator space $\MAX(E)$ over $E$ by letting, for $n \in \posints$ and $x \in M_n(E)$:
\[
  \| x \|_n := \sup \{ \|| x \||_n : \text{$( \|| \cdot \||_m )_{m=1}^\infty$ is a sequence of norms as in Definition \ref{opdef2}} \}.
\]
From this definition it is immediately clear that, for any other operator space $F$, the identity ${\cal B}(E,F) = \CB(\MAX(E),F)$ holds with identical norms. If $\dim E = \infty$, then $\id_E$ does not lie in
$\CB(\MIN(E),\MAX(E))$ (\cite[Theorem 14.3(iii)]{Pau2}): this shows that, generally, there are many different operator space over the same Banach space.
\item Given a Hilbert space $\Hilbert$, we define operator spaces --- called the {\it column\/} and {\it row space\/}, respectively, over $\Hilbert$ --- by letting
\[
  \COL(\Hilbert) := {\cal B}(\comps,\Hilbert) \qquad\text{and}\qquad \ROW(\Hilbert) := {\cal B}(\Hilbert,\comps).
\]
\item Let $E_0$ and $E_1$ be operator spaces such that $(E_0,E_1)$ is a compatible couple of Banach spaces in the sense of interpolation theory (\cite{BL}). Then, for each $n \in \posints$, the couple $M_n(E_0), M_n(E_1))$
is also compatible. For $n \in \posints$ and $\theta \in [0,1]$, we can thus define
\[
  M_n(E_\theta) := ( M_n(E_0), M_n(E_1))_\theta
\]
in the sense of complex interpolation (\cite{Pis1}): this defines an operator space over $E_\theta = (E_0, E_1)_\theta$ (see \cite{Pis1} and \cite{Pis2} for more information).
\end{enumerate}
\par
The most significant advantage the abstract Definition \ref{opdef2} has for us over the concrete Definition \ref{opdef1} is that it allows for the development of a duality theory for operator spaces. Naively, one might think
that the Banach space dual $E^\ast$ of some operator space $E$ can be equipped with an operator space structure through identifying $M_n(E^\ast)$ with $M_n(E)^\ast$; the resulting norms on the spaces $M_n(E^\ast)$, however, will no longer
satisfy (\ref{R1}).
\par
Given an operator space $E$ and $n \in \posints$, we can, for $m \in \posints$, identify $M_m(M_n(E))$ with $M_{mn}(E)$ and thus turn $M_n(E)$ into an operator space. Given two operator spaces $E$ and $F$, we then can use the {\it algebraic\/}
identification
\begin{equation} \label{CBdef}
  M_n(\CB(E,F)) := \CB(E,M_n(F)) \qquad (n \in \posints)
\end{equation}
to equip $\CB(E,F)$ with an operator space structure. We shall use this operator space structure, to turn the Banach space dual of an operator space into an operator space again.
\par
The following is \cite[Theorem 2.10]{Smi} and the starting point for the duality of operator spaces:
\begin{theorem}
Let $E$ be an operator space, and let $n \in \posints$. Then each $T \in {\cal B}(E,M_n)$ is completely bounded such that $\| T \|_\cb = \left\| T^{(n)} \right\|$
\end{theorem} 
\par
Letting $n = 1$, we obtain:
\begin{corollary} \label{smithcor}
Let $E$ be an operator space. Then $E^\ast$ equals $\CB(E,\comps)$ with identical norms.
\end{corollary}
\par
Since $\CB(E,\comps)$ is an operator space by means of (\ref{CBdef}), we can use Corollary \ref{smithcor} to turn $E^\ast$ into an operator space, i.e.\ $M_n(E^\ast) = \CB(E,M_n)$ for $n \in \posints$.
\par
We list a few examples of dual operator spaces:
\begin{examples}
\item Let $E$ be a Banach space. Then we have the canonical completely isometric isomorphisms
\[
  \MIN(E)^\ast = \MAX(E^\ast) \qquad\text{and}\qquad \MAX(E)^\ast = \MIN(E^\ast).
\]
Since the natural operator space structure of a $\cstar$-algebra is always minimal by definition, it follows that the canonical operator space structure of the dual of a commutative $\cstar$-algebra --- or, more generally, of the predual
of a commutative von Neumann algebra --- is always maximal. 
\item Let $\Hilbert$ be a Hilbert space, and let $\overline{\Hilbert}^\ast$ be its conjugate dual space. It is elementary functional analysis that $\Hilbert$ and $\overline{\Hilbert}$ are canonically isometrically isomorphic. For the column
and row spaces over $\Hilbert$, however, we have
\[
  \overline{\COL(\Hilbert)}^\ast = \ROW(\Hilbert) \qquad\text{and}\qquad  \overline{\ROW(\Hilbert)}^\ast = \COL(\Hilbert),
\]
i.e.\ those operator spaces are not self-dual even though the underlying Banach space is a Hilbert space.
\item For each Hilbert space $\Hilbert$, there is a unique operator space $\mathfrak{OH}$ over $\Hilbert$ such that $\overline{\mathfrak OH}^\ast = \mathfrak{OH}$: this operator space was introduced by G.\ Pisier in \cite{Pis1}. It
can we shown (see \cite{Pis1}) that
\[
  \mathfrak{OH} = ( \COL(\Hilbert), \ROW(\Hilbert) )_{\frac{1}{2}} = ( \MIN(\Hilbert), \MAX(\Hilbert) )_{\frac{1}{2}}.
\]
\item Given any measure space $X$, the Banach spaces $L^\infty(X)$ and $L^1(X)$ form a compatible couple such that $(L^\infty(X), L^1(X))_\theta = L^p(X)$ for $p \in [1,\infty]$ and $\theta = \frac{1}{p}$. In \cite{Pis2}, Pisier
used this and the fact that $L^\infty(X)$ --- as a commutative $\cstar$-algebra --- and $L^1(X)$ --- as subspace of the dual of the commutative $\cstar$-algebra $L^\infty(X)$ --- each carry a natural operator space structure to
define an operator space, which we denote by $OL^p(X)$, over $L^p(X)$. For $p,p' \in (1,\infty)$ such that $\frac{1}{p} + \frac{1}{p'} = 1$, we have the duality $OL^p(X)^\ast = OL^{p'}(X)$.
\end{examples}
\par
We conclude this section with the operator space analog of an elementary result on bounded operator between Banach spaces (\cite{Ble1} and \cite[Lemma 1.1]{Ble2}):
\begin{theorem} \label{david}
Let $E$ and $F$ be operator spaces. Then the adjoint of every completely bounded operator from $E$ to $F$ is completely bounded, and the map
\begin{equation} \label{adjoint}
  \CB(E,F) \to \CB(F^\ast, E^\ast), \quad T \mapsto T^\ast
\end{equation}
is a complete isometry.
\end{theorem}
\begin{proof}
It is easy to see (the proof parallels the one in the Banach space situation) that (\ref{adjoint}) is well defined and, in fact, an isometry. To see that (\ref{adjoint}) is even a complete isometry, fix $n \in \posints$,
and note that the following maps and identifications are all isometric:
\begin{eqnarray*}
  M_n(\CB(E,F)) & = & \CB(E,M_n(F)) \\
  & \hookrightarrow & \CB(M_n(F)^\ast, E^\ast), \quad\text{since taking adjoints is an isometry}, \\
  & = & \CB(F^\ast, M_n(E^\ast)) \\
  & = & M_n(\CB(F^\ast, E^\ast)).
\end{eqnarray*}
This proves the claim.
\end{proof}
\par
The book \cite{ER} is the first monograph devoted to the theory of (abstract) operator space: all those assertions for which we did not provide specific references can be found there (along with proofs). Another introduction to
operator spaces by Pisier (\cite{Pis3}) is scheduled to appear soon. V.\ I.\ Paulsen's book \cite{Pau2} has a somewhat different thrust, but also contains the essentials of operator space theory. A little known, hidden gem  --- a strange thing
to say about an article on the world wide web --- is the lexicon style article \cite{Wit}: it focuses on the concepts rather than on detailed technicalities. To get an impression of how operator space theory (or rather the theory of completely
bounded maps between concrete operator spaces) looked like in the pre-Ruan days, see \cite{Pau1}.
\section{Quantized Banach algebras} \label{QBan}
A Banach algebra is an algebra which is also a Banach space such that multiplication is contractive (or merely bounded). To add operator space overtones to that definition, we first have to define what it means for a
bilinear map between operator spaces to be completely bounded.
\par
Let $E_1$ $E_2$, and $F$ be operator spaces, and let $n_1, n_2 \in \posints$. Then the $(n_1,n_2)^{\mathrm{th}}$ amplification of a bilinear map $T \!: E_1 \times E_2 \to F$, denoted by $T^{(n_1, n_2)} \!:
M_{n_1}(E_1) \times M_{n_2}(E_2) \to M_{n_1 n_2}(F)$, is defined as follows: Given $x = \left[ x_{j,k} \right]_{j,k=1}^{n_1} \in M_{n_1}(E_1)$ and $y = \left[ y_{\nu,\mu} \right]_{\nu, \mu=1}^{n_2} \in M_{n_2}(E_2)$, we set
\[
  T^{(n_1, n_2)}(x,y) := \left[ T( x_{j,k}, x_{\nu,\mu}) \right]_{{j,k=1, \ldots, n_1} \atop {\nu, \mu =1, \ldots, n_2}}.
\]
We call $T$ {\it completely bounded\/} if $\|T \|_\cb := \sup_{n_1, n_2 \in \posints} \left\|  T^{(n_1, n_2)} \right\| < \infty$ and {\it completely contractive\/} if $\| T \|_\cb \leq 1$.
\par
We can thus define:
\begin{definition} \label{Qdef}
A {\it quantized Banach algebra\/} is an algebra which is also an operator space such that multiplication is a completely bounded bilinear map.
\end{definition}
\par
We do {\it not\/} demand that multiplication in a quantized Banach algebra be completely contractive (such algebras are called completely contractive Banach algebras in \cite{RuaAG}): all the examples in this section will have
completely contractive multiplication, but we wish to have some more freedom in Section \ref{FTH} below. 
\begin{examples}
\item Let $\A$ be any Banach algebra. Then $\MAX(\A)$ is a quantized Banach algebra.
\item Let $\Hilbert$ be a Hilbert space, and let $\A$ be a closed subalgebra of ${\cal B}(\Hilbert)$. Then $\A$, equipped with its concrete operator space structure, is a quantized Banach algebra: quantized Banach algebras of
this form are called {\it operator algebras\/}. In view of Theorem \ref{Ruan1}, one might jump to the conclusion that every quantized Banach algebra is an operator algebra. This is wrong, however. In fact, the algebras we shall be 
concerned with in Section \ref{QBan2} are rarely operator algebras (we shall discuss this below). An axiomatic description of operator algebras is given in \cite{Ble3} (see also \cite{BRS} for an earlier result in the unital case).
\item Let $E$ be any operator space. Then $\CB(E)$ with the composition of operators as product is a quantized Banach algebra. Such quantized Banach algebras are operator algebras only if $E = \COL(\Hilbert)$ for
some Hilbert space $\Hilbert$, so that canonically 
\[
  \CB(E) = \CB(\COL(\Hilbert)) = {\cal B}(\Hilbert)
\]
as operator spaces (\cite[Theorem 3.4]{Ble3}). 
\end{examples}
\par
As in the Banach space category, there is a projective tensor product of operator spaces, i.e.\ a universal linearizer for bilinear maps between operator spaces (\cite{BP} and \cite{ER1}). 
Following \cite{ER}, we write $\Tensor$ for this tensor product.
Given two operator spaces $E$ and $F$, we may form their operator space tensor product $E \Tensor F$ and their Banach space tensor product $E \tensor^\gamma F$; the universal property of $\tensor^\gamma$ immediately yields a
canonical contraction from $E \tensor^\gamma F$ to $E \Tensor F$. This is about everything that can be said, for general $E$ and $F$, about the relation between $\Tensor$ and $\tensor^\gamma$: the operator space projective tensor product
is {\it not\/} an operator space over the Banach space projective tensor product.
\par
The operator space projective tensor product turns out to be the ``right'' tensor product for the predual of von Neumann algebras in the sense that it enjoys a rather pleasant duality property. 
\par
Given two von Neumann algebras $\M$ and $\N$ acting on Hilbert spaces $\Hilbert$ and $\mathfrak K$, respectively, the {\it von Neumann algebra tensor product\/} $\M \bar{\tensor} \N$ is defined as the von Neumann algebra
acting on the Hilbert space tensor product $\Hilbert \ttensor_2 \mathfrak{K}$ generated by the algebraic tensor product $\M \tensor \N$. The von Neumann algebras $\M$, $\N$, and $\M \bar{\tensor} \N$ each have a (unique) predual
space, $\M_\ast$, $\N_\ast$, and, $(\M \bar{\tensor} \N)_\ast$, respectively. Moreover, we have a canonical map from $\M_\ast \tensor \N_\ast$ to $(\M \bar{\tensor} \N)_\ast$. Since von Neumann algebras have a canonical (concrete)
operator space structure, so have their duals and thus their preduals; hence, we may form the operator space projective tensor product $\M_\ast \Tensor \N_\ast$. As it turns out, this yields a complete description of 
$(\M \bar{\tensor} \N)_\ast$ in terms of $\M_\ast$ and $\N_\ast$ (\cite{ER0}):
\begin{theorem} \label{EdRuan}
Let $\M$ and $\N$ be von Neumann algebras with preduals $\M_\ast$ and $\N_\ast$, respectively. Then we have a canonical, completely isometric isomorphism
\[
  \M_\ast \Tensor \N_\ast \cong (\M \bar{\tensor} \N)_\ast.
\]
\end{theorem}
\par
We shall now exhibit further examples of quantized Banach algebras with the help of Theorem \ref{EdRuan}
\begin{definition} \label{Hopf}
A {\it Hopf--von Neumann algebra\/} is a pair $(\M, \Gamma^\ast)$, where $\M$ is a von Neumann algebra, and
$\Gamma^\ast$ is a {\it co-multiplication\/}: a unital, injective, normal, i.e.\ $w^\ast$-$w^\ast$- continuous, $^\ast$-homomorphism from $\M$ to $\M \bar{\tensor} \M$ which is co-associative, i.e.\ the diagram
\begin{equation} \label{comult}
  \begin{diagram}
  \M                  & \rTo^{\Gamma^\ast}                & \M {\tensor} \M \\                                                                             
  \dTo^{\Gamma^\ast}  &                                   & \dTo_{\Gamma^\ast \tensor \id_\M} \\
  \M \bar{\tensor} \M & \rTo_{\id_\M \tensor \Gamma^\ast} & \M \bar{\tensor} \M \bar{\tensor} \M
  \end{diagram}
\end{equation}
commutes. 
\end{definition}
\par
Let $(\M,\Gamma^\ast)$ be a Hopf--von Neumann algebra. Since $\Gamma^\ast \!: \M \to \M \bar{\tensor} \M$ is $w^\ast$-continuous, it must be the adjoint operator of some $\Gamma \!: (\M_\ast \bar{\tensor} \M_\ast)
\to \M_\ast$. Since $\Gamma^\ast$ as a $^\ast$-ho\-mo\-mor\-phism is a complete contraction, so is $\Gamma$ (as a consequence of Theorem \ref{david}). Invoking Theorem \ref{EdRuan}, we see that $\Gamma$ maps 
$\M_\ast \Tensor \M_\ast$ into $\M_\ast$, thus inducing a completely contractive, bilinear map from $\M_\ast \times \M_\ast$ to $\M_\ast$. The commutativity of the diagram (\ref{comult}) makes sure that
this bilinear map is indeed an associative multiplication on
$\M_\ast$, so that $\M_\ast$ becomes a quantized Banach algebra.
\begin{example}
Let $G$ be a locally compact group. Identifying $L^\infty(G)
\bar{\tensor} L^\infty(G)$ with  $L^\infty(G \times G)$, we define a
co-multiplication $\Gamma^\ast \!: L^\infty(G) \to L^\infty(G)
\bar{\tensor} L^\infty(G)$ by letting
\[
  (\Gamma^\ast \phi)(x,y) := \phi(xy) \qquad (\phi \in L^\infty(G), \, x,y \in G).
\]
Given $f, g \in L^1(G) = L^\infty(G)_\ast$ and $\phi \in L^\infty(G)$,
we obtain that
\begin{eqnarray*}
  \langle \Gamma(f \tensor g), \phi \rangle & = &  \langle f \tensor g, \Gamma^\ast\phi \rangle \\
  & = & \int_{G \times G} f(x)g(y) \phi(xy) \, dx \, dy \\
  & = & \int_{G \times G} f(x) g(x^{-1} y) \phi(y) \, dx \, dy \\
  & = & \int_G (f \ast g)(y) \phi(y) \, dy \\
  & = & \langle f \ast g, \phi \rangle, 
\end{eqnarray*}
i.e.\ the multiplication on $L^1(G)$ induced by $\Gamma^\ast$ is
nothing but the convolution product (\ref{conv}). Since the von
Neumann algebra $L^\infty(G)$ is commutative, the underlying operator
space of this quantized Banach algebra is $\MAX(L^1(G))$.
\end{example}
\par
This example shows that a quantized Banach algebra, even one arising
as the predual of a Hopf--von Neumann algebra need not be an operator
algebra: If $L^1(G)$ were an operator algebra, it would, in particular, be a closed
subalgebra of an Arens regular Banach algebra and therefore be Arens
regular itself; this, however, is possible only if $G$ is finite (\cite{You}).
\par
Next, we shall see that not only $L^1(G)$, but also $A(G)$ and $B(G)$
are quantized Banach algebras in a canonical manner.
\begin{example}
Let $G$ be a locally compact group, and $\lambda$ be the left regular
representation of $G$ on $L^2(G)$. The {\it group von Neumann
  algebra\/} $\VN(G)$ of $G$ is defined as $\lambda(G)''$. The {\it
  fundamental operator\/} $W \in {\cal B}(L^2(G))$, defined through
\[
  (W\xi)(x,y) := \xi(x,xy) \qquad (\xi \in L^2(G), \, x,y \in G),
\]
is easily seen to be unitary. Letting
\[
  \Gamma^\ast T := W^\ast(T \tensor \id_{L^2(G)})W \qquad (T \in \VN(G)),
\]
we obtain a co-multiplication $\Gamma^\ast \!:\VN(G) \to \VN(G)
\bar{\tensor} \VN(G)$, thus turning $\VN(G)_\ast$ into a
quantized Banach algebra. The predual of $\VN(G)$, however, is nothing
but $A(G)$ (\cite{Eym}): Given $T \in \VN(G)$ and $f \in A(G)$ as in
(\ref{Fourier3}), the duality is implemented via $\langle f, T \rangle
:= \langle T\xi, \eta \rangle$. Since 
\[
  \Gamma^\ast\lambda(x) = \lambda(x) \tensor \lambda(x) \qquad (x \in G),
\]
the multiplication on $A(G)$ induced by $\Gamma^\ast$ is 
pointwise multiplication.
\end{example}
\par 
For non-discrete or for infinite, amenable $G$, the Fourier algebra
$A(G)$ fails to be Arens regular (\cite{For}) and thus cannot be an
operator algebra. The operator space structure of $A(G)$ is further
investigated in \cite{FW}. In the same paper, it is also observed that
the canonical operator space structure of $A(G)$ is different from
$\MAX(A(G))$ unless $G$ is abelian.
\par
Given two locally compact groups $G$ and $H$, it is easy to see that
$\VN(G) \bar{\tensor} \VN(H) \cong \VN(G \times H)$ in a canonical
manner. In view of Theorem \ref{EdRuan}, we thus obtain the following
extremely useful tensor identity:
\begin{corollary} \label{AGcor}
Let $G$ and $H$ be locally compact groups. Then we have a canonical
completely isometric isomorphism
\[
  A(G) \Tensor A(H) \cong A(G \times H).
\]
\end{corollary}
\par
The corresponding tensor identity for $\tensor^\gamma$ is false in
general: In fact, V.\ Losert (\cite{Los}) showed that $A(G)
\tensor^\gamma A(H) \cong A(G \times H)$ holds isomorphically if and
only if $G$ or $H$ has an abelian subgroup of finite index; this
isomorphism is an isometry if and only if $G$ or $H$ is
abelian. 
\par
Concluding this section, we turn to $B(G)$:
\begin{example}
Let $G$ be a locally compact group, and let $\pi_u$ be its universal
representation. Let $\wstar(G) := \pi_u(G)''$, i.e.\ the second dual
of $\cstar(G)$. The representation $\pi_u$ has the following
universal property: For any representation $\pi$ of $G$ on a Hilbert space, there is a unique normal $^\ast$-homomorphism
$\rho \!: \wstar(G) \to \pi(G)''$ such that $\pi = \rho \circ \pi_u$. Applying this universal property to the representation
\[
  G \to \wstar(G) \bar{\tensor} \wstar(G), \quad x \mapsto \pi_u(x) \tensor \pi_u(x)
\]
yields a co-multiplication $\Gamma^\ast \!: \wstar(G) \to \wstar(G)
\bar{\tensor} \wstar(G)$. Hence, $\cstar(G)^\ast$ is a quantized
Banach algebra. From the definitions of $\cstar(G)$ and $B(G)$,
however, it is clear that $B(G)$ and $\cstar(G)^\ast$ can be
canonically identified and that the multiplication on $B(G)$ induced
by $\Gamma^\ast$ is pointwise multiplication.
\end{example}
\par
The two operator space structures on $A(G)$ --- the one it has as the
predual space of $\VN(G)$ and the one it inherits from $B(G)$ --- are
identical.
\par
The Hopf--von Neumann algebras $L^\infty(G)$ and $\VN(G)$ have
additional structure making them (more or less the only) examples of
{\it Kac algebras\/} (see \cite{ES} for the precise definition). Kac
algebras have a duality theory which extends the well known Pontryagin
duality for locally compact, abelian groups (see again \cite{ES} for
details). For applications of operator spaces to the study of abstract
Kac algebras, see \cite{RuaKac}, \cite{KR1}, \cite{KR2}, and \cite{RX}, for instance.
\section{Homological properties of Fourier and Fou\-rier--Stiel\-tjes algebras} \label{QBan2}
Like the notion of a Banach algebra, the concept of a Banach module
also translates painlessly into the quantized world:
\begin{definition}
A bimodule $E$ over a quantized Banach algebra $\A$ is called a {\it quantized Banach $\A$-bimodule\/} if it is also an operator space such that the module actions
\[
  \A \times E \to E, \quad 
   (a,x) \mapsto \left\{ \begin{array}{c} a \cdot x, \\ x \cdot a
   \end{array} \right.
\] 
are completely bounded.
\end{definition}
\par
It is routinely checked that $E^\ast$ with its dual operator space
structure and the dual module actions is again a quantized Banach 
$\A$-bimodule.
\par
In analogy with Definition \ref{bamdef}, we may thus define (following
\cite{RuaAG}):
\begin{definition} \label{opamdef}
A quantized Banach algebra $\A$ is called {\it operator amenable\/} if, for each quantized Banach $\A$-bimodule $E$, every completely bounded derivation $D \!: \A \to E^\ast$ is inner.
\end{definition}
\par
The theory of operator amenable, quantized Banach algebra unfolds
parallel to the usual theory of amenable Banach algebras: The
intrinsic characterization Proposition \ref{diagprop} holds true in the
quantized category as well (with $\tensor^\gamma$ replaced by
$\Tensor$) as does the collection of hereditary properties Theorem
\ref{bamher} (all bounded maps have to replaced by completely bounded
ones) --- the proofs carry over almost verbatim. In fact, the whole
theory of amenable Banach algebras can be viewed as a subset of the
theory of operator amenable, quantized Banach algebras: A Banach
algebra $\A$ is amenable if and only if the quantized Banach algebra
$\MAX(\A)$ is operator amenable.
\par
The following theorem due to Z.-J.\ Ruan (\cite[Theorem 3.6]{RuaAG})
put the concept of operator amenability on the mathematical map:
\begin{theorem} \label{zhongjin}
Let $G$ be a locally compact group. Then the following are equivalent:
\begin{items}
\item $G$ is amenable.
\item $A(G)$ is operator amenable.
\end{items}
\end{theorem}
\begin{proof}
(i) $\Longrightarrow$ (ii): postponed.
\par
(ii) $\Longrightarrow$ (i): Operator amenable, quantized Banach
algebras always have bounded approximate identities: this is proven in
exactly the same manner as Proposition \ref{bambai}. Hence,
Theorem \ref{lep} yields the amenability of $G$.
\end{proof}
\par
Comparing Theorem \ref{zhongjin} and \ref{Barry2}, we see that operator
amenability establishes a much more satisfactory correspondence
between the structure of $A(G)$ and properties of $G$ than does
amenability in the sense of Definition \ref{bamdef}.
\par
Of course, like amenability, other homological concepts such as weak
amenability, biprojectivity, and biflatness can be provided with
operator space overtones. How is this is done is straightforward and
left to the reader.
\par
In analogy with Proposition \ref{biproj}, a quantized Banach algebra
$\A$ is operator biprojective if and only if the multiplication map
$\Gamma \!: \A \Tensor \A \to \A$ has a completely bounded right
inverse which is also an $\A$-bimodule homomorphism.
\par
In view of Theorem \ref{sasha0} and Pontryagin duality, one might expect that $A(G)$ is operator biprojective if and only if $G$ is discrete. This is indeed the case as was independently shown
by O.\ Yu.\ Aristov (\cite{Ari}) and P.\ J.\ Wood (\cite{Woo}). 
\par
We first record a lemma:
\begin{lemma} \label{lem1}
Let $\A$ be a commutative, operator biprojective, quantized Banach algebra with character space $\Phi_\A$. Then $\Phi_\A$ is discrete.
\end{lemma}
\par
The proof of the classical counterpart (e.g., \cite[Corollary 2.8.42]{Dal}) carries over with the obvious modifications.
\begin{theorem} \label{oleg1}
The following are equivalent for a locally compact group $G$:
\begin{items}
\item $G$ is discrete.
\item $A(G)$ is operator biprojective.
\end{items}
\end{theorem}
\begin{proof}
(i) $\Longrightarrow$ (ii): Let $\chi_\Gamma$ be the indicator function of the {\it diagonal subgroup\/}
\[
  G_\Gamma := \{ (x,x) : x \in G \}.
\]
Then $\chi_\Gamma$ is positive definite and, since $G$ is discrete, lies in $B(G \times G)$. Define
\[
  \rho \!: A(G) \to B(G \times G), \quad f \mapsto (f \tensor 1)\chi_\Gamma.
\]
Then $\rho$ is completely bounded, and it is easy so see that it attains its values in $A(G \times G)$, i.e.\ in $A(G) \Tensor A(G)$ by Corollary \ref{AGcor}. It is routinely
checked that $\rho$ is a bimodule homomorphism and a right inverse of $\Gamma$.
\par
(ii) $\Longrightarrow$ (i) is clear by Lemma \ref{lem1}.
\end{proof}
\par
In analogy with the situation for Banach algebras, operator biflatness is weaker then both operator amenability and operator biprojectivity. In our discussion of the group algebra $L^1(G)$, we noted that
biflatness for such algebras is the same as amenability. This is not the case for $A(G)$: the Fourier algebra is operator biprojective --- and thus, in particular, operator biflat --- whenever $G$ is discrete, so that,
for example, $A(\free_2)$ is operator biflat, but not operator amenable. Since the class of groups $G$ for which $A(G)$ is operator biflat includes all amenable and all discrete groups, it is a tempting conjecture that
$A(G)$ is operator biflat for {\it every\/} locally compact group.
\par
The following proposition is from \cite{ARS}:
\begin{proposition} \label{ARS1}
Let $G$ be a locally compact group, and suppose that there is a bounded net $( \mathbf{f}_\alpha )_\alpha$ in $B(G \times G)$ with the following properties:
\begin{alphitems}
\item $\lim_\alpha f (\mathbf{f}_\alpha |_{G_\Gamma}) = f$ for all $f \in A(G)$;
\item $\lim_\alpha \mathbf{g} \mathbf{f}_\alpha = 0$ for all $\mathbf{g} \in \ker \Gamma$.
\end{alphitems}
Then $A(G)$ is operator biflat.
\end{proposition}
\begin{proof}
Let $\mathbf{F} \in B(G \times G)^{\ast\ast}$ be a $w^\ast$-accumulation point of $( \mathbf{f}_\alpha )_\alpha$. Then
\[
  \rho \!: A(G) \to B(G \times G)^{\ast\ast}, \quad f \mapsto f \cdot \mathbf{F}
\]
can be shown to attain its values in $A(G \times G)^{\ast\ast}$. It is routinely verified that $\rho$ is an $A(G)$-bimodule homomorphism as required by (the quantized counterpart of) Definition \ref{biflat}.
\end{proof}
\par
This leaves us with the question of whether a net as required by Proposition \ref{ARS1} always exists.
\par
Recall that a locally compact group is a $[\SIN]$-group if $L^1(G)$ has a bounded approximate identity, $( e_\alpha )_\alpha$ say, belonging to its center i.e.\ satisfying $\delta_x \ast e_\alpha = e_\alpha \ast \delta_x$
for all indices $\alpha$ and for all $x \in G$ with $\delta_x$ denoting the point mass at $x$. Every discrete group is trivially a $[\SIN]$-group. 
\par
We define:
\begin{definition}
A locally compact group $G$ is called a {\it quasi-$[\SIN]$-group\/} if $L^1(G)$ has a bounded approximate identity $( e_\alpha )_\alpha$ such that
\begin{equation} \label{qsinlim}
  \delta_x \ast e_\alpha - e_\alpha \ast \delta_x \to 0 \qquad (x \in G).
\end{equation}
\end{definition}
\par
All $[\SIN]$-groups are trivially quasi-$[\SIN]$-groups, but so are all amenable groups; a connected group is even quasi-$[\SIN]$ if and only if it is amenable (\cite{LR}). 
\par
The following theorem is from \cite{RX}, but the proof, which avoids the Kac algebra machinery used in \cite{RX}, is from \cite{ARS}, and even yields a slightly stronger result (\cite[Theorem 2.4]{ARS}):
\begin{theorem} \label{ARS2}
Let $G$ be a quasi-$[\SIN]$-group. Then $A(G)$ is operator biflat.
\end{theorem}
\begin{proof}
By \cite{LR} and \cite{Sto}, we can find a bounded approximate identity $( e_\alpha )_\alpha$ for $L^1(G)$ such that:
\begin{itemize}
\item $e_\alpha \geq 0$ and $\| e_\alpha \|_1 = 1$ for each index $\alpha$;
\item the limit (\ref{qsinlim}) is uniform on compact subsets of $G$;
\item for each neighborhood $U$ of $e$ there is an index $\alpha_U$ such that $\supp \, e_\alpha \subset U$ for all $\alpha \succ \alpha_U$.
\end{itemize}
Let $\lambda, \rho \!: G \to {\cal B}(L^2(G))$ be the left and right
regular representation, respectively, of $G$. Letting
\[
  \mathbf{f}_\alpha(x,y) := \left\langle \lambda(x) \rho(y) e_\alpha^\frac{1}{2}, e_\alpha^\frac{1}{2} \right\rangle \qquad (x,y \in G),
\]
we obtain a net $( \mathbf{f}_\alpha )_\alpha$ as required by Proposition \ref{ARS1}, so that $A(G)$ is operator biflat.
\end{proof}
\par
With Theorem \ref{ARS2} proven, we can now complete the proof of Theorem \ref{zhongjin}: If $G$ is amenable, it is a quasi-$[\SIN]$-group by \cite{LR}, so that $A(G)$ is operator biflat. Since $A(G)$ has a bounded approximate identity
by Theorem \ref{lep} due to the amenability of $G$, it follows from (the quantized counterpart of) Theorem \ref{sasha} that $A(G)$ is operator amenable.
\par
Nevertheless, the question of whether $A(G)$ is always operator biflat remains open: For $G = \operatorname{SL}(3,\comps)$, it can be shown that a net as in Proposition \ref{ARS1} does not exists (\cite[Theorem 4.5]{ARS}) even though
this does not rule out that $A(G)$ is operator biflat for some other reason.
\par
It is well known that biflatness implies weak amenability (\cite[Proposition 2.8.62]{Dal} or \cite[Theorem 5.3.13]{LoA}): In analogy, operator biflatness implies operator weak amenability. As it turns out, $A(G)$ is at least always
operator amenable (\cite{Spr}):
\begin{theorem} \label{nico}
Let $G$ be a locally compact group. Then $A(G)$ is operator weakly amenable.
\end{theorem}
\begin{proof}
We only deal with the compact case, so that $A(G)$ has an identity; the general case is dealt with by adjoining an identity, but this somewhat complicates the argument.
\par
In view of \cite{Groe} or \cite{Run0}, it is sufficient that $(\ker \Gamma)^2$, i.e.\ the linear span of all product of elements in $\ker \Gamma$, is dense in $\ker \Gamma$.
Identifying $A(G) \Tensor A(G)$, with $A(G \times G)$ the multiplication operator $\Gamma$ becomes the restriction to $G_\Gamma$, i.e.\
\[
  \ker \Gamma = \{ \mathbf{f} \in A(G \times G) : \mathbf{f} |_{G_\Gamma} = 0 \}.
\] 
Since $G_\Gamma$ is a subgroup of $G \times G$, it is a set of synthesis for $A(G \times G)$ by \cite[Theorem 2]{Her}. Since the closure of $(\ker \Gamma)^2$ and $\ker \Gamma$ clearly have the same hull, they must be equal.
\end{proof}
\par
Theorems \ref{zhongjin}, \ref{oleg1}, and \ref{nico} are rather satisfactory in the sense that they provide beautiful insights into the homological nature of the Fourier algebra in terms of the underlying group. Nevertheless,
the questions of whether $A(G)$ is, e.g.\ amenable in the sense of Definition \ref{bamdef} or weakly amenable, are equally legitimate and shouldn't be dismissed as uninteresting just because the probable answers are likely not to be
as pleasant as in the quantized setting.
\par
In \cite{FKLS}, the operator space structure of $A(G)$ is used to answer a question that {\it a priori\/} seems to have nothing to do with operator spaces. The authors attempt to determine which closed ideals of $A(G)$
have a bounded approximate identity. For amenable $G$, they characterize those ideals as being precisely those of the form
\[
  I(E) = \{ f \in A(G) : f |_E = 0 \},
\]
where $E$ is a set in the closed coset ring of $G$: this considerably generalizes the abelian case dealt with in \cite{LRW}.
\par
An important step in their argument is the following result:
\begin{proposition}
Let $G$ be an amenable, locally compact group, and let $H$ be a closed subgroup of $G$. Then $I(H)$ has a bounded approximate identity.
\end{proposition}
\begin{proof}
By Theorem \ref{zhongjin}, $A(G)$ is operator amenable. The annihilator $I(H)^\perp$ in $\VN(G)$ of $I(H)$ is easily seen to be the $w^\ast$-closed linear span of the set $\{ \lambda(x) : x \in H \}$, which can be identified
with $\VN(H)$. Since $G$ is amenable, so is $H$, which, in turn, implies that the von Neumann algebra $\VN(H)$ is injective (see \cite[Chapter 6]{LoA}). In particular, there is a norm one projection ${\cal E} \!:
\VN(G) \to \VN(H)$. The projection $\cal E$ can then be shown to be completely bounded, so that $I(H)^\perp$ is completely complemented in $\VN(G)$. From the (quantized counterpart of) Theorem \ref{bamher}(iii), it then follows
that $I(H)$ is operator amenable and thus has a bounded approximate identity.
\end{proof}
\par
Using the results from \cite{FKLS} and some more operator space theory, the question for which locally compact groups $G$ precisely $A(G)$ is amenable can be settled (\cite{FR}):
\begin{theorem} \label{Runde}
Let $G$ be a locally compact group. Then the following are equivalent:
\begin{items}
\item $G$ has an abelian subgroup of finite index.
\item $A(G)$ is amenable.
\end{items}
\end{theorem}
\begin{proof}
We have already seen that (i) $\Longrightarrow$ (ii) holds (Theorem \ref{LLWthm}). 
\par
We only give a brief outline for the proof of the converse. If $A(G)$ is amenable, then it can be shown that the ideal
\[
  I = \{ f \in A( G \times G) : \text{$f(x,x^{-1}) = 0$ for $x \in G$} \}
\]
has a bounded approximate identity, so that, by \cite{FKLS}, the set $\{ (x,x^{-1}) : x \in G \}$ must lie in the coset ring of $G \times G$.
By an old result from \cite{Rud}, this means that
\[
  G \to G, \quad x \mapsto x^{-1}
\]
is what is called piecewise affine, and a recent result due to M.\ Ilie (\cite{Ili}), then entails that
\[
  A(G) \to A(G), \quad f \mapsto \check{f}
\]
is completely bounded. This possible only if $G$ is as in Theorem \ref{Runde}(i).
\end{proof}
\par
As for the result from \cite{FKLS}, the {\it statement\/} of Theorem \ref{Runde} makes no reference to the operator space structure of $A(G)$ whereas its proof relies on it.
\par
As far as the weak amenability of $A(G)$ is concerned, there is a plausible, but still open conjecture: $A(G)$ is weakly amenable if and only if the principal component of $G$ is abelian. In \cite{FR},
Forrest and the author have shown the sufficiency of this condition.
\par
So far, the only quantized Banach algebra we have dealt with in this section was $A(G)$. We now turn briefly to $B(G)$.
\par
In view of Theorem \ref{DGHthm} and the abelian case, the canonical conjecture is that $B(G)$ is operator amenable if and only if $G$ is compact. So, far the best result in the difficult direction of this conjecture is
given in \cite{RS}:
\begin{theorem} \label{nicoandi}
The following are equivalent for a locally compact group:
\begin{items}
\item $G$ is compact.
\item $B(G)$ is operator $C$-amenable for some $C < 5$, i.e.\ $B(G)$ has an approximate diagonal in $B(G) \Tensor B(G)$ bounded by some $C < 5$.
\end{items}
\end{theorem}
\par
The direction (i) $\Longrightarrow$ (ii) is fairly easy: Since $G$ is compact, $B(G)$ equals $A(G)$ and thus is operator amenable by Theorem \ref{zhongjin}. An inspection of the proof shows that the amenability of $G$ does, in fact,
already imply the operator $1$-amenability of $A(G)$. The proof of the hard direction (ii) $\Longrightarrow$ (i) makes use of a decomposition of $B(G)$ that can be interpreted as an analogue of the decompsition (\ref{dirsum}) of
$M(G)$ into its discrete and continuous parts.
\par
Somewhat surprisingly, the question for which $G$ the Fourier--Stieltjes algebra is amenable as a Banach algebra has a complete answer (\cite{FR}); it is a fairly easy corollary of Theorem \ref{Runde}:
\begin{corollary} 
The following are equivalent for a locally compact group $G$:
\begin{items}
\item $G$ has a compact, abelian subgroup of finite index.
\item $B(G)$ is amenable.
\end{items}
\end{corollary}
\begin{proof}
(i) $\Longrightarrow$ (ii): In this case, $G$ itself is compact, so that $B(G) = A(G)$ must be amenable by Theorem \ref{LLWthm}.
\par
(ii) $\Longrightarrow$ (i): Since $A(G)$ is a complemented ideal in $B(G)$, it is also amenable. Hence, by Theorem \ref{Runde}, $G$ has an abelian subgroup $H$ of finite index. Without loss of generality let $H$ be closed
and thus open, so that the restriction map from $B(G)$ to $B(H)$ is surjective. Hence, $B(H)$ is also amenable. Since $B(H) \cong M(\hat{H})$ via the Fourier--Stieltjes transform, the measure algebra $M(\hat{H})$ must be amenable. 
Amenability of $M(\hat{H})$, however, forces $\hat{H}$ to be discrete (Theorem \ref{DGHthm}). Hence, $H$ must be compact.
\end{proof}
\section{Fig\`a-Talamanca--Herz algebras} \label{FTH}
Let $G$ be a locally compact group, let $p \in (1,\infty)$, and let $p' \in (1,\infty)$ be dual to $p$, i.e.\ $\frac{1}{p} + \frac{1}{p'} = 1$. The {\it Fig\`a-Talamanca--Herz algebra\/} $A_p(G)$ consists of those functions $f \!: G \to \comps$
such that there are sequences $( \xi_n )_{n=1}^\infty$ in $L^p(G)$ and $( \eta_n )_{n=1}^\infty$ in $L^{p'}(G)$ with
\begin{equation} \label{Apeq1}
  \sum_{n=1}^\infty \| \xi_n \|_{L^p(G)} \| \eta_n \|_{L^{p'}(G)} < \infty
\end{equation}
and
\begin{equation} \label{Apeq2}
  f = \sum_{n=1}^\infty \xi_n \ast \check{\eta}_n.
\end{equation}
The norm of $f \in A_p(G)$ is defined as the infimum over all sums (\ref{Apeq1}) such that (\ref{Apeq2}) holds. It is clear that $A_p(G)$ is a quotient space of $L^p(G) \tensor^\gamma L^{p'}(G)$  that embeds contractively into ${\cal C}_0(G)$. 
For $p =2$, we obtain the Fourier algebra $A(G)$.
\par
The following is \cite[Theorem B]{Her0}:
\begin{theorem} \label{carl1}
Let $G$ be a locally compact group, and let $p,q \in (1,\infty)$ be such that $p \leq q \leq 2$ or $2 \leq q \leq p$. Then pointwise multiplication induces a contraction from $A_p(G) \tensor^\gamma A_q(G)$ to $A_p(G)$.
\end{theorem}
\par
For $p = q$, this implies:
\begin{corollary}
Let $G$ be a locally compact group, and let $p \in (1,\infty)$. Then $A_p(G)$ is a Banach algebra under pointwise multiplication.
\end{corollary}
\par
If $G$ is amenable, $A_p(G)$ has an approximate identity that can be chosen to be bounded by $1$: the proof of Theorem \ref{lep} from \cite{Lep} carries over to this more general setting. Hence, we obtain (\cite[Theorem C]{Her0}):
\begin{corollary}
Let $G$ be an amenable locally compact group, and let $p,q \in (1,\infty)$ be such that $p \leq q \leq 2$ or $2 \leq q \leq p$. Then $A_p(G)$ is contained in $A_p(G)$ such that the inclusion map is a contraction.
\end{corollary}
\par
In view of Theorems \ref{bamher}(i) and \ref{LLWthm}, this yields:
\begin{corollary}
Let $G$ be a locally compact group with an abelian subgroup of finite index. Then $A_p(G)$ is amenable for all $p \in (1,\infty)$.
\end{corollary}
\par
Let $\lambda_{p'} \!: G \to {\cal B}(L^{p'}(G))$ be the regular left representation of $G$ on $L^{p'}(G)$. Via integration, $\lambda_{p'}$ extends to a representation of $L^1(G)$ on $L^{p'}(G)$. 
The algebra of {\it $p'$-pseudomeasures\/} $\PM_{p'}(G)$ is defined as the $w^\ast$-closure of $\lambda_{p'}(L^1(G))$ in ${\cal B}(L^{p'}(G))$. There is a canonical duality $\PM_{p'}(G) \cong A_p(G)^\ast$ via
\[
  \langle \xi \ast \check{\eta}, T \rangle := \langle T\eta, \xi \rangle \qquad (\xi \in L^{p'}(G), \, \eta \in L^p(G), \, T \in \PM_{p'}(G)). 
\]
For $p =2$, this is just the usual duality between $A(G)$ and $\VN(G)$.
\par
If we want to use operator space techniques to deal with $A_p(G)$, we are faced with a problem right at the beginning: Beside $\MAX(A_p(G))$, which is uninteresting as an operator space, there seems to be no operator space over
$A_p(G)$ turning it into a quantized Banach algebra. The operator space structure of $A(G)$ stems from its duality with $\VN(G)$, whose operator space structure is a concrete one, arising from $\VN(G) \subset {\cal B}(L^2(G))$.
Since ${\cal B}(L^2(G)) = \CB(\COL(L^2(G)))$, we may want to attempt to construct a column operator space over $L^{p'}(G)$ and use the duality between $A_p(G)$ and $\PM_{p'}(G) \subset \CB(\COL(L^{p'}(G)))$ to equip $A_p(G)$ with an
operator space structure.
\par
At the first glance, column and row spaces make no sense for Banach spaces other than Hilbert spaces. The following characterization of column and row Hilbert spaces, however, due to B.\ Matthes (\cite{Math}), indicates a way to
circumvent this difficulty:
\begin{theorem} \label{maththm}
Let $\Hilbert$ be an operator space whose underlying Banach space is a Hilbert space. Then the following are equivalent:
\begin{items}
\item $\Hilbert = \COL(\Hilbert)$.
\item $M_{n,1}(\Hilbert) = M_{n,1}(\MAX(\Hilbert))$ and $M_{1,n}(\Hilbert) = M_{1,n}(\MIN(\Hilbert))$ holds for all $n \in \posints$, i.e.\ $\Hilbert$ is maximal on the columns and minimal on the rows.
\end{items}
\end{theorem}
\par
A similar characterization holds for $\ROW(\Hilbert)$.
\par
In his doctoral thesis under G.\ Wittstock's supervision (\cite{Lam}), A.\ Lambert uses Theorem \ref{maththm} to define column and row spaces over {\it arbitrary\/} Banach spaces. His crucial idea is to introduce an intermediate category
between Banach and oprator spaces, the so-called {\it operator sequence spaces\/}:
\begin{definition} \label{ossdef}
An {\it operator sequence space\/} is a linear space $E$ with a complete norm $\| \cdot \|_n$ on $E^n$ for each $n \in \posints$ such that
\begin{eqnarray*} 
  \| (x,0) \|_{n+m} & = & \| x \|_n \qquad (n,m \in \posints, \, x \in E^n), \\ 
  \| (x,y) \|^2_{n+m} & \leq & \| x \|_n^2 + \| y \|_m^2 \qquad (n,m \in \posints, \, x \in E^n, \, y \in E^m)
\end{eqnarray*}
and
\[
  \| \alpha x  \|_m \leq \| \alpha \| \| x \|_n \qquad (n,m \in \posints, \, x \in E^n, \, \alpha \in M_{m,n}).
\]
\end{definition}
\par
In analogy with the completely bounded maps, one can define appropriate morphisms for operator sequence spaces called {\it sequentially bounded maps\/} and denoted by $\SB(E,F)$ for two operator sequence spaces $E$ and $F$. 
As for operator spaces here is a duality for operator sequence spaces, which allows to equip the Banach space dual of an operator sequence space again with an operator sequence space structure 
(see \cite{Lam} or \cite{LNR} for this and more). 
\begin{examples}
\item Let $E$ be an operator space. Identifying, $E^n$ with $M_{n,1}(E)$ for each $n \in \posints$, i.e.\ taking the columns of $E$, we obtain an operator sequence space which we denote by $C(E)$.
\item Given a Banach space $E$, we may identity $E^n$ with ${\cal B}(\ell^2_n,E)$ and obtain the {\it minimal operator sequence space\/} $\min(E)$ over $E$. For each other operator sequence space $F$, we have
${\cal B}(F,E) = \SB(F,\min(E))$ with identical norms, which justifies the adjective ``minimal''.
\item Let $E$ be any Banach space. The {\it maximal operator sequence space\/} $\max(E)$ over $E$ is defined as follows: For $n \in \posints$ and $x \in E^n$, define
\[
  \| x \|_n := \inf \{ \| \alpha \| \| y \|_{\ell^2_m(E)} : m \in \posints, \, \alpha \in M_{n,m}, \, y \in E^m, \, x = \alpha y \}.
\]
This operator sequence space has the property that ${\cal B}(E,F) = \SB(\max(E),F)$ --- with identical norms --- for any other operator sequence space $F$.
\end{examples}
\par
The operator sequence spaces $\min$ and $\max$ are dual to one another, i.e.\
\[
  \min(E)^\ast = \max(E^\ast) \qquad\text{and}\qquad \max(E)^\ast = \min(E^\ast)
\]
for each Banach space $E$. The proof parallels the one for the corresponding assertion about $\MIN$ and $\MAX$. 
\par
Even though the basics of operator sequence spaces very much parallel the corresponding results about operator spaces, the category of operator sequences spaces sometimes displays phenomena putting it
closer to Banach spaces:
\begin{enumerate}
\item As we have already noted (\cite[Theorem 14.3]{Pau2}), the identity on a Banach space $E$ is not completely bounded from $\MIN(E)$ to $\MAX(E)$ if $E$ is infinite-dimensional. If $\A$ is a $\cstar$-algebra, however, then
$\id_\A$ lies in $\SB(\min(\A),\max(\A))$ if and only if $\A$ is subhomogeneous (\cite[Satz 2.2.25]{Lam}).
\item The principle of local reflexivity, which is a cornerstone of the local theory of Banach spaces, but fails to have an analog for general operator spaces, still works in the category of operator sequence spaces (\cite[Satz 1.3.26]{Lam}).
\end{enumerate}
\par
The first of our examples shows that the columns of an operator space always form an operator sequence space. As we shall now see, this is the only way an operator sequence space can occur, but possibly in more than one fashion:
\begin{definition} \label{Mindef}
Let $E$ be an operator sequence space. Then the {\it minimal operator space\/} $\Min(E)$ over $E$ is defined by letting $M_n(\Min(E)) := {\cal B}( \ell^2_n, E^n)$.
\end{definition}
\par
The adjective ``minimal'' is again justified in the usual way: If $F$ is any other operator space, then we have an isometric identity $\SB(C(F),E) = \CB(F,\Min(E))$. It follows that $\Min(\min(E)) = \MIN(E)$ for each Banach space $E$.
\par
\begin{definition} \label{Maxdef}
Let $E$ be an operator sequence space. Then the {\it maximal operator space\/} $\Max(E)$ over $E$ is defined by letting, for $x \in M_n(E)$,
\[
  \| x \|_n := \inf \{ \| \alpha \| \| \beta \| : x = \alpha \, \diag(v_1, \ldots, v_k) \beta \},
\]
where the infimum is taken over all $k,l \in \posints$, $\alpha \in M_{n,kl}$, $\beta \in M_{k,n}$, and $v_1, \ldots, v_k$ in the closed unit ball of $E^l$.
\end{definition}
\par
Given any other operator space $F$, we then have the isometric identity $\SB(E,C(F)) = \CB(\Max(E),F)$, so that the adjective ``maximal'' makes sense. Moreover, $\Max(\max(E)) = \MAX(E)$ holds for any Banach space $E$, and $\Min$ and $\Max$
are dual to one another, i.e.\
\[
  \Min(E)^\ast = \Max(E^\ast) \qquad\text{and}\qquad \Max(E)^\ast = \Min(E^\ast)
\]
holds for each operator sequence space $E$.
\par
With Theorem \ref{maththm}, we can now quote the definition --- from \cite{Lam} --- of $\COL(E)$ and $\ROW(E)$ where $E$ is an arbitrary Banach space:
\begin{definition} \label{coldef}
Let $E$ be a Banach space.
\begin{alphitems}
\item The {\it column space\/} over $E$ is defined as $\COL(E) := \Min(\max(E))$.
\item The {\it row space\/} over $E$ is defined as $\ROW(E) := \Max(\min(E))$.
\end{alphitems}
\end{definition}
\par
It follows from Theorem \ref{maththm} that, for a Hilbert space, Definition \ref{coldef} yields the usual column and row spaces.
\begin{proposition} \label{colprop}
Let $E$ be a Banach space. Then $\COL(E)$ and $\ROW(E)$ are operator spaces such that
\[
  {\cal B}(E) = \CB(\COL(E)) = \CB(\ROW(E)) 
\]
with identical norms and
\[
  \COL(E)^\ast = \ROW(E^\ast) \qquad\text{and}\qquad \ROW(E)^\ast = \COL(E^\ast).
\]
\end{proposition}
\par
In view of the properties of $\min$, $\max$, $\Min$, and $\Max$ and their various dualities, the verification of Proposition \ref{colprop} is fairly straightforward.
\par
We now return to the Fig\`a-Talamanca--Herz algebras.
\par
Given a locally compact group $G$ and $p,p' \in (1,\infty)$ dual to one another, we use the inclusion
\[
  \PM_{p'}(G) \subset {\cal B}(L^{p'}(G)) = \CB(\COL(L^{p'}(G)))
\]
to define an operator space over $\PM_{p'}(G)$. Via the duality $A_p(G)^\ast = \PM_{p'}(G)$, we then obtain an operator space structure on $A_p(G)$.
\par
For this operator space structure, we obtain the following (\cite[Theorem 6.4]{LNR}):
\begin{theorem}
Let $G$ be a locally compact group, and let $p,q \in (1,\infty)$ be such that $p \leq q \leq 2$ or $2 \leq q \leq p$. Then pointwise multiplication induces a completely bounded map from $A_p(G) \Tensor A_q(G)$ to $A_p(G)$.
\end{theorem}
\par
We do not know if the map from  $A_p(G) \Tensor A_q(G)$ to $A_p(G)$ is even a complete contraction; an upper bound for its $\cb$-norm is given in \cite{LNR}.
\par
The following corollaries are immediate:
\begin{corollary}
Let $G$ be a locally compact group, and let $p \in (1,\infty)$. Then $A_p(G)$ is a quantized Banach algebra under pointwise multiplication.
\end{corollary}
\begin{corollary} \label{incl}
Let $G$ be an amenable locally compact group, and let $p,q \in (1,\infty)$ be such that $p \leq q \leq 2$ or $2 \leq q \leq p$. Then $A_p(G)$ is contained in $A_p(G)$ such that the inclusion map is completely bounded.
\end{corollary}
\par
We can now extend Theorem \ref{zhongjin} to Fig\`a-Talamanca--Herz algebras:
\begin{theorem} \label{opam}
The following are equivalent for a locally compact group $G$:
\begin{items}
\item $G$ is amenable.
\item $A(G)$ is operator amenable.
\item $A_p(G)$ is operator amenable for each $p \in (1,\infty)$.
\item There is $p \in (1,\infty)$ such that $A_p(G)$ is operator amenable.
\end{items}
\end{theorem}
\begin{proof}
(i) $\Longleftrightarrow$ (ii) is Theorem \ref{zhongjin}.
\par
(ii) $\Longrightarrow$ (iii): Let $p \in (1,\infty)$. If $A(G)$ is operator amenable, then $G$ is amenable, so that $A(G) \subset A_p(G)$ completely boundedly by Corollary \ref{incl}. Since this inclusion has dense range,
the hereditary properties of operator amenability imply the operator amenability of $A_p(G)$.
\par
(iii) $\Longrightarrow$ (iv) is trivial.
\par
(iv) $\Longrightarrow$ (i): Let $p \in (1,\infty)$ be such that $A_p(G)$ is operator amenable; in particular, $A_p(G)$ then has a bounded approximate identity. Since Theorem \ref{lep} holds for $A_p(G)$ as well, we conclude that
$G$ is amenable.
\end{proof}
\par
It is only natural to ask which of the results proved for the quantized Banach algebra $A(G)$ in the previous section carry over to Fig\`a-Talamanca--Herz algebra. In the proofs of Theorems \ref{oleg1} and \ref{nico}, the tensor identity
from Corollary \ref{AGcor} plays a pivotal r\^ole. We therefore conclude our survey with the following open problem:
\begin{problem}
Let $G$ and $H$ be locally compact groups, and let $p \in (1,\infty)$ be arbitrary. Do we have a canonical, completely bounded (but not necessarily completely isometric) isomorphism
\[
  A_p(G) \Tensor A_p(H) \cong A_p(G \times H)
\]
as in the case $p =2$? (By a completely bounded isomorphism we mean an isomorphism which is completely bounded with a completely bounded inverse.) 
\end{problem}
\dated
\vfill
\begin{tabbing}
{\it Author's address\/}: \= Department of Mathematical and Statistical Sciences \\
\> University of Alberta \\
\> Edmonton, Alberta \\
\> Canada T6G 2G1 \\[\medskipamount]
{\it E-mail\/}: \> {\tt vrunde@ualberta.ca} \\[\medskipamount]
{\it URL\/}: \> {\tt http://www.math.ualberta.ca/$^\sim$runde/}
\end{tabbing}           
\end{document}

%% file: format.tex
\renewcommand{\baselinestretch}{1.2}
\newcommand{\dated}{\mbox{} \hfill {\small [{\tt \today}]}}

%% file: mathdefs.tex
\usepackage{amsmath,amssymb,amscd}
\newcommand{\pf}[1]{\trivlist \item[\hskip\labelsep\it #1\ ]}
\newcommand{\varpf}[1]{\trivlist \item[\hskip\labelsep\sc #1:]}
\newcommand{\qedbox}{$\rlap{$\sqcap$}\sqcup$}
\newcommand{\qed}{\qquad \qedbox \endtrivlist}
\newcommand{\varqed}{\hfill \rule{0.6em}{0.6em} \endtrivlist}
\newenvironment{proof}{\pf{Proof}}{\qed}

\newenvironment{example}{\pf{Example}}{\endtrivlist}
\newenvironment{examples}{\pf{Examples} 
   \begin{enumerate}}{\end{enumerate} \endtrivlist}
\newenvironment{items}{
  \begin{enumerate} 
                    
  }{\end{enumerate}}
\newenvironment{alphitems}{
  \begin{enumerate} 
                    
  }{\end{enumerate}}
\newenvironment{keywords}{\noindent\small {\it Keywords\/}:}{\vskip 4pt}
\newenvironment{classification}{\noindent\small 2000 {\it Mathematics Subject
Classification\/}:}{\vskip 12pt}

\newcommand{\comps}{{\mathbb C}}
\newcommand{\reals}{{\mathbb R}}
\newcommand{\ints}{{\mathbb Z}}
\newcommand{\posints}{{\mathbb N}}

\newcommand{\torus}{{\mathbb T}}
\newcommand{\free}{{\mathbb F}}

\newcommand{\tensor}{\otimes}
\newcommand{\ttensor}{\tilde{\otimes}}
\newcommand{\Tensor}{\hat{\otimes}}

\newcommand{\cstar}{{C^\ast}}
\newcommand{\wstar}{{W^\ast}}
\newcommand{\re}{{\operatorname{Re}}}
\newcommand{\im}{{\operatorname{Im}}}
\newcommand{\id}{{\mathrm{id}}}

\newcommand{\diag}{{\operatorname{diag}}}

\newcommand{\A}{{\mathfrak A}}
\newcommand{\B}{{\mathfrak B}}
\newcommand{\Hilbert}{{\mathfrak H}}
\newcommand{\M}{{\mathfrak M}}
\newcommand{\N}{{\mathfrak N}}

\newcommand{\VN}{\operatorname{VN}}

\newcommand{\SIN}{\operatorname{SIN}}

\newcommand{\supp}{{\operatorname{supp}}}

%% file: thmdefs.tex
\newtheorem{theorem}{Theorem}[section]
\newtheorem{lemma}[theorem]{Lemma}
\newtheorem{corollary}[theorem]{Corollary}
\newtheorem{proposition}[theorem]{Proposition}
\newtheorem{df}[theorem]{Definition}
\newenvironment{definition}{\begin{df} \rm}{\end{df}}